\documentclass{amsproc}

\usepackage{amsfonts,amsmath,amssymb,epsf, epsfig}
\usepackage{mathrsfs}
\usepackage{graphics} 
\usepackage{ifpdf}
\ifpdf
\usepackage{epstopdf}
\usepackage{hyperref}
\else
\usepackage[hypertex]{hyperref}
\fi
\usepackage{xypic}
\usepackage{epigraph}
\usepackage{psfrag}
\usepackage{hhline}

\DeclareGraphicsExtensions{.ps}

\newtheorem{thm}{Theorem}[section]

\newtheorem{defn}[thm]{Definition}
\newtheorem{prop}[thm]{Proposition}

\newtheorem{rema}[thm]{Remark}

\newtheorem{conj}[thm]{Conjecture}

\numberwithin{equation}{section}

\newcommand\void[1]       {}

\newcommand\be            {\begin{equation}}
\newcommand\bea           {\begin{eqnarray}}
\newcommand\bearll        {\begin{array}{ll}\displaystyle}

\newcommand\ee            {\end{equation}}
\newcommand{\eea}         {\end{eqnarray}}
\newcommand\eear          {\end{array}}

\newcommand\bnu           {\begin{enumerate}}
\newcommand\enu           {\end{enumerate}}

\newcommand\etb           {&\!\! \displaystyle}

\newcommand{\nn}          {\nonumber \\}

\newcommand\cl           {\mathrm{cl}}
\newcommand\op         {\mathrm{op}}
\newcommand\co         {\text{\rm cl-op}}
\newcommand\pt          {\mathrm{pt}}
\newcommand\CcC     {\Cc_V^{(2)}}

\renewcommand\cir         {\,{\circ}\,}

\newcommand\eps           {\varepsilon}
\newcommand\Hom           {\mathrm{Hom}}
\newcommand\id            {{\rm id}}

\newcommand\one           {{\bf1}}
\newcommand\oti           {\,{\otimes}\,}
\newcommand\ti            {\,{\times}\,}
\newcommand\wt          {\mathrm{wt}}
\newcommand\lwt         {\mathrm{wt}^L}
\newcommand\rwt        {\mathrm{wt}^R}

\newcommand\Cb            {\mathbb{C}}
\newcommand\Hb           {\mathbb{H}}
\newcommand\Nb            {\mathbb{N}}

\newcommand\Rb            {\mathbb{R}}
\newcommand\Zb            {\mathbb{Z}}

\newcommand\Cc            {\mathcal{C}}
\newcommand\DC            {\mathcal{D}}
\newcommand\Ec            {\mathcal{E}}
\newcommand\Fc           {\mathcal{F}}
\newcommand\Gc          {\mathcal{G}}

\newcommand\LC            {\mathcal{L}}
\newcommand\calM           {\mathcal{M}}
\newcommand\Nc            {\mathcal{N}}

\newcommand\Y             {\mathcal{Y}}

\begin{document}

\title{Conformal field theory and a new geometry}

\author{Liang Kong}

\address{Institute for Advanced Study, Tsinghua University, Beijing 100084 China}

\email{kong.fan.liang@gmail.com}

\thanks{The author is supported in part by the Basic Research Young Scholars Program of Tsinghua University, Tsinghua University independent research Grant No. 20101081762 and by NSFC Grant No. 20101301479. }

\subjclass{Primary: 81T40; Secondary: 17B69, 18D10, 18D50, 18R10}
\date{July 14, 2011}

\keywords{Conformal field theory, D-branes, vertex operator algebra, stringy algebraic geometry}

\begin{abstract} 
This paper is a review of open-closed rational conformal field theory (CFT) via the theory of vertex operator algebras (VOAs), together with a proposal of a new geometry based on CFTs and D-branes. We will start with an outline of the idea of the new geometry, followed by some philosophical background behind this vision. Then we will review a working definition of CFT slightly modified from Segal's original definition and explain how VOA emerges from it naturally. Next, using the representation theory of rational VOAs, we will discuss a classification result of open-closed rational CFTs, from which some basic properties of a rational CFT, such as the Holographic Principle, can be derived. They will also serve as supporting evidences for the vision of a new geometry. In the end, we briefly discuss the connection between our vision of a new geometry and other topics. 
\end{abstract}

\maketitle

\setcounter{footnote}{0}
\def\thefootnote{\arabic{footnote}}

\bigskip

\small{
\begin{quote}
\hspace{4.03cm}{\it To see a world in a grain of sand,  

\hspace{3.4cm}And a heaven in a wild flower,  

\hspace{3.4cm}Hold infinity in the palm of your hand, 

\hspace{3.4cm}And eternity in an hour.} \\

\hspace{3.4cm} --- William Blake ``Auguries of Innocence" 
\end{quote}
}

\section{Introduction and summary}   \label{sec:intro}

This paper has grown out of many talks I have given in 2009 and 2010 on a new geometry based on conformal field theory (CFT). Although the idea is rooted in many works in both physics and mathematics, as far as I know, it has never been clearly stated and emphasized. It is not surprising that this sounds quite shocking and alien to some of my audiences, even to myself when I first realized it in May 2007. So in this paper, I will try to cover some philosophical aspects of this new geometry which are impossible to cover in my talks. Besides this vision of a new geometry, this paper is pretty much a review of open-closed CFT from the point of view of vertex operator algebra (VOA). To avoid any confusion, CFTs in this paper are all 2-dimensional. Moreover, a CFT in Section \ref{sec:intro} and \ref{sec:geometry} means a CFT with or without supersymmetry, but for later sections, we will only discuss CFT without supersymmetry because the mathematical foundations of open-closed supersymmetric CFT are still lacking. 

\medskip
String theory is proposed to be a theory of quantum gravity. By Einstein, gravity is nothing but geometry. Therefore, it is equally good to say that string theory is a proposal for a theory of {\em quantum geometry}. One of the aims of string theory is to answer how geometry emerges. All classical notions in geometry, such as points, lines, surfaces, dimensions and metrics, should all be derived concepts. They should arise as certain derived quantities, moduli or invariants from the abstract mathematical structures of string theory. By giving up all familiar geometric notions, string theory challenges us to find a radical new point of view of geometry. 

It is impossible to know how to write a new geometry from nothing. Any theoretical speculation needs hints from experimental data. One of the testing grounds or laboratories of this new quantum geometry is provided by the so-called non-linear sigma models. These models connect mathematical structures of string theory directly to notions in classical geometry. The outcomes of this laboratory are many ``phenomenological theories" manifesting themselves as new topics in geometry. To name a few:  mirror symmetry, quantum cohomology, Gromov-Witten theory, topological string, elliptic cohomology, etc. What is the secret structure behind of all these clues?

Perhaps we can ask an easier question. What is the new message from string theory that is so distinguishable that it alone can explain why we missed so many new insights before the advent of string theory? There are perhaps many answers to this question. One answer is that string theory seems to emphasize a point of view from loop space which has richer structures than the original manifold. If we want to take this message seriously, we would immediately ask the question: what are the structures on loop space? Again, there are perhaps many answers to this question. One answer suggested by string theory is that the free loop space, or rather the space of functions on the free loop space, has a natural algebraic structure which is called closed conformal field theory. One should take the term ``loop space" and ``CFT" in a flexible way.  They can have different meanings in different contexts. For example, loop space can be defined in both the topological context and the (derived-)algebro-geometric context.

\medskip
CFT is a formidable subject and unfamiliar to most mathematicians. In spite of many works on this subject, its full nature is still very mysterious even to experts. Ironically, within the frameworks of known definitions of CFT, we cannot even claim to have a single nontrivial example! Nevertheless, many ingredients of CFT such as vertex operator algebra, modular functor, modular tensor category, etc. are explicitly constructed and relatively well understood. Moreover, many variants of CFT or VOA such as string topology \cite{string-top}, topological CFT (TCFT) \cite{costello}, homological CFT \cite{godin}, ..., etc. have been constructed. Very often, for convenience, we will pretend to talk about a CFT while we are only talking about its substructures such as a VOA or its variants such as a TCFT or string topology.

Inspired by the suggestion that the free loop space is a closed CFT, a few groups of mathematicians initiated programs to realize this picture. We will list some constructions by these mathematicians. The list is not complete due to my limited knowledge.
\begin{itemize}
\item  Malikov, Schechtman and Vaintrob constructed the Chiral de Rham complex \cite{msv} which is a sheaf of vertex operator algebras (VOA) on a smooth manifold. It can be viewed as a shadow of certain structure on formal loop space. 

\item Kapranov and Vasserot \cite{kv1} constructed an algebro-geometric version of free loop space and showed that it has a natural structure of factorization monoid, which can be viewed as a non-linear version of factorization algebra that contains VOA as local data. 

\item Chas and Sullivan constructed the closed string topology \cite{string-top}, which is an algebraic structure on the homology of the free loop space. A closed string topology can be viewed as a homological closed conformal field theory \cite{godin}. 

\item Ben-Zvi, Francis and Nadler \cite{bfn} showed that the stable symmetric monoidal $\infty$-category of quasi-coherent sheaves on the loop space of a perfect stack has the structure of a 2-dimensional topological field theory (TFT) \cite{lurie-tft} which is a categorified analogue of  a TCFT \cite{costello}. 

\end{itemize}
By all these works, we can be more confident about identifying a closed CFT as an algebraic model for a free loop space, at least for our philosophical discussion. More importantly, as a theory of quantum geometry, string theory demands us to forget about the loop space in the first place because it can only be a derived concept. Instead, we shall take a closed CFT as initial data to rebuild the entire geometry. 

One convenient approach to see the connection between CFT and geometry is to work only with certain substructures of nonlinear sigma models. For example, substructures such as superconformal algebras, chiral rings, A-branes, B-branes, partition functions, etc., are directly linked to mirror symmetry, quantum cohomology, Fukaya category, bounded derived category of coherent sheaves, elliptic genus, respectively. This approach has been very successful so far. But the disadvantage of it is that we might lose the global picture (if there is any) and perhaps some deeper revelations.

\medskip
In recent years, there have been some developments on the mathematical foundations of open-closed rational CFTs without supersymmetry\footnote{I believe that supersymmetric CFTs share similar properties with CFTs without supersymmetry. We will not distinguish them in our philosophical discussion in this section.}  via the representation theory of VOAs \cite{osvoa}\cite{ffa}\cite{ko-ffa}\cite{ocfa}\cite{cardy-cond}, and another successful and nearly complementary program led by Fuchs, Runkel and Schweigert \cite{fs-cat}\cite{tft1}\cite{tft3}\cite{tft4}\cite{tft5} in which an open-closed rational CFT is studied as a holographic boundary of a 3-dimensional topological field theory. These two approaches provide compatible results \cite{tft1}\cite{cardy-alg} and can be unified \cite{cardy-alg-2}. Together, they give a rather complete\footnote{High genus theory for rational CFTs is still not available. But they are uniquely determined by genus-0 theory and there are good reasons to believe that they indeed exist.} picture of rational open-closed CFTs. Perhaps, it is dangerous to use the word ``complete" because we are always subject to certain axiom systems. But I believe that the system is rich enough to reveal many global features of CFT and to recover or link to other missing aspects. As we will review in later sections, these developments, in particular, give the following two results:  
\bnu
 
\item A closed CFT is indeed a stringy generalization of a commutative ring. As we will show in Section \ref{sec:ffa}, a closed CFT is a full field algebra satisfying a stringy commutativity condition, and also a commutative associative algebra in a braided tensor category \cite{ko-ffa} (see Theorem \ref{thm:genus-0-cft}). 

\begin{figure}
  \begin{picture}(200,220)
     \put(-50,0){\scalebox{.6}{\includegraphics{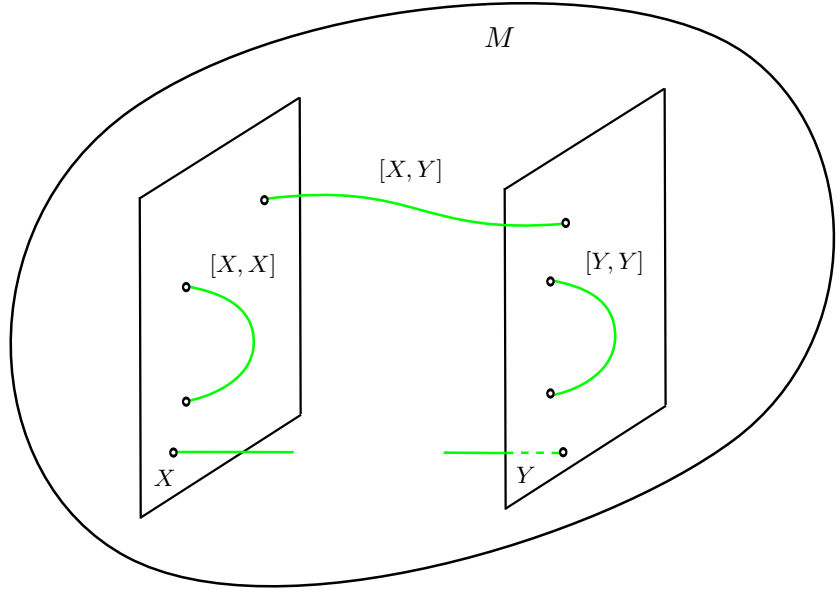}}}
   \put(90, 156){$ [X,Y] $}
   \put(26, 120){$ [X, X] $}
   \put(168, 121){$ [Y,Y] $}
   \put(130, 205){\large $ M $}
    \put(5, 39){$X$}
   \put(139, 40){ $Y$ }
  \end{picture}
  \caption{ {\small Geometric intuition of the category of D-branes associated to a target manifold $M$: 
  a closed CFT $A_\cl$ models the free loop space; a boundary condition, i.e. a chiral $A_\cl$-module $X$, might be able to model the space of $[0,1)$-paths starting from a submanifold associated to $X$ (also denoted by $X$); an open CFT $[X,X]$ models the space of paths from $X$ to $X$; and the $[Y,Y]$-$[X,X]$-bimodule $[X,Y]$ models the space of paths from $X$ to $Y$. }} 
 \label{fig:d-brane}
\end{figure}

\item For a closed CFT $A_\cl$, its D-branes can be parametrized by modules over an open CFT compatible with $A_\cl$. It will become clear later that these modules are automatically {\it chiral modules}\footnote{The word ``chiral" refers to the fact that it is a module over only one copy of the Virasoro algebra (see Remark \ref{def:chiral-module}).} over $A_\cl$ (see Remark \ref{rema:cl-cft-mod=dbrane}, \ref{rema:non-cardy}). 
\enu
I believe that this relation between a  closed CFT and its D-branes is a stringy analogue of that between a commutative ring and its prime ideals. Moreover, it is well known to physicists that D-branes, as boundary conditions for open strings, indeed behave like generalized points or submanifolds\footnote{Usually, such a submanifold is equipped with additional structures, such as a bundle over it and a connection.}, e.g. $X$ and $Y$ in a target manifold $M$ in Figure \ref{fig:d-brane}. This fact has been used by physicists to probe the geometry of the target manifold. What is suggested by above results is a surprisingly simple picture: CFT provides an immediate stringy generalization of classical algebraic geometry!

\medskip
In order to get a more precise picture of this analogy, we have to take a closer look at the nature of the ``set" of D-branes. Because the operator product expansion (OPE) of fields on a D-brane provides a stringy algebra structure, a module $X$ over an open CFT $A$ does not contain the information of fields on a D-brane in general. It is only a boundary condition literally. This fact becomes rather clear in Fuchs-Runkel-Schweigert's state-sum construction of rational CFTs via 3-d TFT techniques \cite{tft1}, in which boundary edges of a triangulation of a surface are labeled by modules over 
an open CFT $A$ which is a separable symmetric Frobenius algebra in a modular tensor category $\Cc$ (see Section \ref{sec:cat-ffa}, \ref{sec:examples}). As we will show later in this case, the closed CFT is given by the center $Z(A)$ of $A$. An $A$-module is automatically a chiral module over the closed CFT. Moreover, the category of $A$-modules is an indecomposable component of the category of chiral modules over $Z(A)$ (see Theorem \ref{thm:open-exists}). Very often, {\it boundary conditions} and {\it D-branes} are used synonymously in the literature. In this work, we will distinguish them. For us, a boundary condition is just an $A$-module $X$ for an open CFT $A$ and a chiral module over the closed CFT. Since the geometric intuition behind a closed CFT $A_\cl$ is given by the free loop space $LM$, I believe that the geometric intuition behind a chiral $A_\cl$-modules $X$ is the space of $[0,1)$-paths starting from a submanifold associated to $X$\footnote{A $[0,1)$-path starting from a submanifold $X$ of $M$ is a continuous map $\gamma: [0,1) \to  M$ such that $\gamma(0)\in X$.} (see Figure \ref{fig:d-brane}). 
I believe that one can find a precise meaning of it in the framework of string topology. 
We further assume that the open CFT $A$ can be chosen to be commutative in a stringy sense\footnote{I believe that such stringy commutative open CFTs exist in most non-linear sigma models.} so that the category of boundary conditions $\Cc_A$ is a monoidal category and we can talk about internal homs\footnote{For irrational CFTs, a tensor product may not exist. Nevertheless, due to the existence of the OPE, at least certain weak version of a tensor product should exist so that we can still talk about internal homs.}. In Section \ref{sec:examples}, we will give an example of such situation which is called Cardy case.

A D-brane is a physical (or dynamical) object, the physical content of which lies in its interaction with other D-branes (including itself) and the bulk theory. From the intuition of the loop space, 
the information of the ``self-interaction" of a D-brane  is contained in the fusion of all open strings with two ends ending on it. For example, when the boundary condition is just a point, all the open strings ending on it form a based loop space. The algebraic model for the fusion of open strings is given by an open CFT compatible with the given closed CFT. By ``compatible" we mean that together they form an open-closed CFT (see Section \ref{sec:def-oc-cft} for a precise definition). As we will show later on, in a rational CFT, for a boundary condition $X$, the corresponding open CFT is given by the internal hom $[X, X]$ in $\Cc_A$. In this case, the closed CFT is also given by the center of $[X,X]$. For irrational CFTs, due to the existence of OPE, it is reasonable to believe that such internal-hom type of constructions still work. Geometrically, it is perhaps also true that one can view the space of paths from a submanifold $X$ to $X$ as some kind of ``internal endo-hom" of the space of $[0,1)$-paths starting from $X$. More general ``interaction" between two different D-branes $X, Y$ is given by the space of paths between them. Its algebraic model is given by a $[Y,Y]$-$[X,X]$-bimodule $[X,Y]$ which is also an internal hom in $\Cc_A$. In the end, we obtain a category of boundary conditions but enriched by internal homs\footnote{It is a standard construction from a module $\calM$ over a tensor category $\Cc$ to obtain a category $\calM^\Cc$ enriched by internal homs in $\Cc$.}.  We will call this enriched category {\it the category of D-branes}\footnote{Notice that the category of D-brane is not much different from a single open CFT because the direct sum (or a coend in general) of the hom spaces of the category of D-branes can gives again an open CFT.}. The geometric intuition behind the category of D-branes is depicted in Figure \ref{fig:d-brane}. A D-brane can then be defined by the functor $[X, -]$. Sometimes, it is also convenient, by only looking at its relation to the bulk theory, to define a D-brane by a pair $(X, [X, X])$ or just the open CFT $[X,X]$. In the rest of this work, we will ignore the ambiguity in the meaning of a D-brane, and leave readers to figure out which one is appropriate in the context. For a given boundary condition $X$, the map $Y \mapsto [X, Y]$ defines a functor from the category of boundary conditions to the category of right $[X, X]$-modules. This functor is an equivalence for all $X$ in a rational CFT. Nothing is known for irrational CFTs, but we assume that the equivalence holds at least for some $X$. That is why we often see in the literature that the category of boundary conditions is defined to be the category of modules over an open CFT. 

As we have mentioned, objects in the category of D-branes can be viewed as stringy primary ideals of the closed CFT which is a stringy analogue of a commutative ring. A chiral module over a closed CFT is a module over only one copy of the Virasoro algebra (or a superconformal algebra in a superconformal field theory) which can be roughly viewed as the Laplacian (or the Dirac operator) on loop space or path spaces. It means that the hom spaces of the category of D-branes are physical spectra. Therefore, the category of D-branes is naturally a combination of the algebro-geometric spectrum (boundary conditions) and physical spectra (internal homs). We will take it to be a stringy analogue of the spectrum of commutative rings\footnote{Of course, it makes little difference to use the category of boundary conditions instead to define the spectrum because the relation between these two categories is clear and standard.} in usual algebraic geometry. This suggests to lift the usual algebraic description of geometry as familiar from the field of algebraic geometry to a {\it stringy algebraic geometry} (SAG), which is summarized in the following table:
\medskip
\begin{center}
   \begin{tabular}{| p{5.4cm} | p{5.3cm} |}
      \hline
      Classical algebraic geometry  & Stringy algebraic geometry  \\  \hhline{|=|=|}
      A commutative ring $R$ & a closed CFT $A_{\mathrm{cl}}$  \\  \hline
      $\text{Spec}(R)$ = the set of prime ideals of $R$ & $\text{Spec}(A_{\mathrm{cl}})$ = the category of D-branes    
       \\  \hline
   \end{tabular}
\end{center}
\medskip
If we agree that such a new geometry indeed exists, the SAG will have the following new features that are different from classical algebraic geometry:  
\bnu

\item CFT emphasizes a bordism point of view towards the notion of space instead of the usual sheaf-theoretical point of view.  It regards a space as a network of subspaces linked by paths instead of a union of open charts. In particular, the spectrum is a category instead of a set. 

\item If we take a closed CFT loosely as an algebraic model for free loop space and an open CFT as an model for path space with a given boundary condition, we see immediately that this stringy geometry has the so-called holographic phenomenon. For example, a D-brane associated to a point is just the based loop space which contains the information of the entire space. It is very different from the usual sheaf-theoretical geometry in which a single chart does not contain any global information. This aspect is somewhat similar to that of noncommutative geometry \cite{connes}, where the noncommutativity encodes the information of gluing data\footnote{I want to thank Matilde Marcolli for clarifying this point.}.  
 
\item It includes spectral geometry as an ingredient. The information of ``distance" and ``intersection" between two D-branres $X$ and $Y$ should be all encoded in the internal hom $[X, Y]$ which is the spectrum of Laplacian or Dirac operator. Be aware that the usual notion of distance only makes sense in the low energy limit. At high energy, different points are highly entangled via large loops or paths. Therefore, we expect a radical new point of view of metric in SAG. 

\enu
We will further discuss the naturalness of these categorical, holographic and spectral nature of SAG in Section \ref{sec:geometry} in order to support our vision. We now summarize these properties of SAG and their classical counterparts in the following table: 
\medskip
\begin{center}
\small{   \begin{tabular}{| p{2.5cm} | p{2.5cm} | p{2.5cm} | p{2.5cm} |}
      \hline
      Classical geometry  & Classical algebra & Stringy geometry & Stringy algebra \\  \hhline{|=|=|=|=|}
      Affine scheme $M$ & comm. ring $R$ & loop space $LM$ & closed CFT  \\  \hline
      points $x, y\in M$ & prime ideals of $R$  & submanifolds of $X, Y\subset M$ & boundary conditions $X, Y$. \\ \hline
      & 
      & the path space between $X$ and $Y$ & internal hom $[X,Y]$  \\  \hline
        & & Laplacian (Dirac) operators on $LM$ or path spaces & Virasoro algebra (superconformal algebra) \\  \hline
   \end{tabular}
}
\end{center}
\medskip
From the content and the properties of SAG, it is fair to say that SAG is a natural combination of algebraic geometry and spectral geometry.

It is tempting to enrich the ``stringy geometry" column by including ``higher paths" and higher bordisms between higher paths all the way to infinity to make it an $\infty$-category \cite{lurie}. We stop at level 2 because of the very nature of string theory and 2-dimensional CFT. There is no higher data on the CFT side to encode the higher bordisms. If one insists to consider higher paths beyond level 2, one need consider higher-dimensional QFTs, which provide higher categorical spectra. It is unclear how a CFT can be extended to a higher dimensional theory. But this might be due to our limited understanding of the situation. Recently, Morrison and Walker proposed a scheme to extend lower dimensional QFTs to higher dimensional QFTs \cite{walker}. I hope to come back to this issue in the future. 

The classical notion of manifold should emerge in the low energy regime as moduli spaces of certain substructures or invariants of the spectrum of SAG. For example, the target manifold can be emergent as the moduli space of certain 0-dimensional objects in the category of D-branes \cite{aspinwall}; the time evolution might come from the moduli space of certain defects (see Section \ref{sec:basic-SAG}). Another way to let a manifold emerge from CFTs is to take the so-called large volume limit. It was studied at a mathematical level by Roggenkamp-Wendland \cite{rw1,rw2} and Kontsevich-Soibelman \cite{konts-soib} (see also Soibelman's contribution \cite{soibelman} in this book). In \cite{konts-soib,soibelman}, it was shown that collapsing a family of unitary CFTs leads to Riemannian manifolds (possible singular) with non-negative Ricci curvature.

\medskip
As many differences between SAG and the notion of scheme are noticeable, one can question the necessity of the notion of spectrum of this new geometry. Even if we take it for granted, one can still question the naturalness of such a definition of spectrum\footnote{Another category plays an important role in CFT is the bicategory with open CFTs as objects and the bimodules over two open CFTs as 1-cells and bimodule maps as 2-cells. Equivalently, one can define an object by the category of modules (also called boundary conditions) over an open CFT, and define a hom category by the category of functors. This bicategory appears often in the study of extended TFTs. In an $n$-dimensional TFT, it becomes the $n$-category of $(n-1)$-categories of boundary conditions. Comparing to the category of D-branes, the essential new things in this bicategory is the inclusion of all bimodules. There are physical motivations to include more general bimodules. They are called defects in physics. Some of them encode the information of dualities \cite{ffrs-duality} of the bulk theory (see Section \ref{sec:defect-duality}).}.
Without any real progress in examples of this new geometry, it is hard to tell what is the natural thing to do at the current stage. An ideal attitude is to keep an open mind. I hope that such vagueness will not be viewed as an evidence against its existence. In Section \ref{sec:geometry}, we will try to analyze the issue from the historical, philosophical and physical point of view in order to provide some support for our vision of SAG.

\medskip
To the best of my knowledge, above simple picture has never been clearly stated in the literature. It is, however, certainly not new to physicists, and has been actively pursued by many string theorists. In the years before my sudden enlightenment in May, 2007 during a talk given by Tom Bridgeland in Max-Planck Institute for Mathematics at Bonn, I was very much influenced by Connes' noncommutative geometry \cite{connes}, Douglas' D-geometry \cite{douglas}, Aspinwall's stringy geometry as the moduli space of D0-branes \cite{aspinwall}, the geometric interpretation of D-branes in non-linear Sigma models \cite{fffs}\cite{bdlr}\cite{mms} and Kapranov-Vasserot's program on infinite dimensional algebraic geometry \cite{kv1}\cite{kapranov2}\cite{kv2}. 

The idea that the rich structure of a CFT can be reduced to produce a spectral triple \cite{connes} goes back to Fr\"{o}hlich and Gawedzki \cite{fg}, later developed by Roggenkamp and Wendland \cite{rw1}, and more recently it was shown that a proper completion of a super-VOA leads to a spectral triple \cite{chkl}\footnote{Recently, Carpi, Hillier, Kawahigashi, Longo and Xu \cite{chklx} connect subfactor theory and noncommutative geometry through the conformal net approach towards superconformal field theory. D-branes are not mentioned in their approach. I hope that a program of new geometry parallel to SAG can be developed in the conformal-net approach to CFT.}. But by doing so, some stringy features of CFT are lost. I believe that there should be an intrinsic stringy geometry associated to CFT.  Such a picture was suggested by Kontsevich and Soibelman. They proposed to look at a CFT as a stringy generalization of a spectral triple \cite{konts-soib} (see Soibelman's contribution \cite{soibelman} in this book). For this reason, SAG can very well be called {\it 2-spectral geometry}\footnote{This name was suggested by Urs Schreiber. I would also like to thank him for pointing out the works \cite{konts-soib}\cite{soibelman} to me after I sent him the first draft of this paper.}. So our proposal on SAG is not new. But Kontsevich-Soibelman's proposal seems to emphasize different aspects of the same picture. 

SAG proposed here has a strong commutative flavor because a closed CFT is stringy commutative. Perhaps that is why the low energy effective theory of closed string theory can describe classical gravity, which is a commutative geometry and refuses to be quantized in any naive way\footnote{Quantizing gravity is not necessary the right question to ask in the first place \cite{hu}.}. The physics on a D-brane, as governed by an open CFT, is noncommutative in general even in the stringy sense. Thus we expect its low energy effective theory to be quantum in nature. 

During my preparation of this paper, I came across an interesting article ``{\em What is a brane?}" by Gregory Moore \cite{Mo2}, in which he wrote: 
\begin{quote}
{\small However, a common theme in the study of D-branes has been the idea that in fact, the (string) field theory on the brane is the primary concept, whereas the spacetime itself is a secondary, derived, concept. This notion has been given some degree of precision in the so-called Matrix theory formulation of M-theory. A rough analogy of what physicists expect may be described in the context of purely topological branes, where the field theory on a brane is described in terms of a noncommutative Frobenius algebra, and the
ÒspacetimeÓ in which it propagates is derived from the Hochschild cohomology of that algebra. These
ideas might ultimately lead to a profound revision of the way we regard spacetime.}
\end{quote}
As we will show in later sections that a closed rational CFT is a commutative Frobenius algebra in a braided tensor category and a rational D-brane (an open CFT) is a noncommutative Frobenius algebra in a tensor category. Moreover, the closed CFT very often appears as the center of the open CFT. Meanwhile, Hochschild cohomology is nothing but a notion of derived center. Therefore, the proposal of SAG can be viewed as a relatively detailed explanation of Moore's prediction in the framework of CFT from my own point of view.

\medskip
There are other attempts in proposing new geometries associated to D-branes. One influential approach is the so-called derived algebraic geometry \cite{dag}, which was also influenced at its early stage by the intuition from CFT through Kontsevich. In a nutshell, derived algebraic geometry replaces a commutative ring by the Hochschild cohomology $\mathrm{HH}^\bullet(A)$ of an associative algebra $A$, which carries a structure of topological conformal field theory \cite{costello}. Therefore, it can be viewed as a parallel program to our stringy algebraic geometry. In general, one can replace the associative algebra by an $A_\infty$-algebra, or a $\otimes_\infty$-category, or an $E_n$-algebra, etc. DAG is a very rich theory, but it contains no information on the metric. Another algebraic geometry approach towards D-branes was taken by Liu and Yau \cite{liuyau1,liuyau2}, they proposed an Azumaya non-commutative geometry formulation of D-branes. This approach is sheaf-theorectical and largely motivated by formulating new moduli problems. The relation between D-branes and schemes was also discussed by Gomez and Sharpe \cite{gomez-sharpe}.

\medskip
Perhaps it is wise to remind readers again that it is still premature to say what this new geometry really is because of the lack of examples. Nevertheless, it is important to initiate the discussion on a larger scale about the general philosophy of this topic as such philosophical inquiries have already generated new questions and inspired new works including some of my own works on CFT. As we will shown in Section \ref{sec:outlooks}, it is perhaps also interesting to explore the possible application of such philosophy to other subjects such as string topology, geometric Langlands correspondence, etc. The SAG proposed here is a step towards such goal. Similar to Connes' noncommutative geometry, working with examples is crucial in this program. We will pursue it in our future publications.

\medskip
The layout of this paper: in Section 2, we will discuss some historical, philosophical and physical backgrounds of SAG; in Section 3, we will review the mathematical foundations of closed CFTs; in Section 4, we will review the mathematical foundations of open-closed CFTs, and give a classification of open-closed CFTs over a rational VOA; in Section 5, we will discuss a few properties (relevant to SAG) of open-closed rational CFTs; in Section 6, we will briefly discuss the connection between SAG and other topics.

\medskip
Convention of notations: $\Nb, \Zb, \Zb_+, \Rb, \Rb_+, \Cb$ denote the set of natural numbers, integers, positive integers, real numbers, positive real numbers, complex numbers, respectively. 
We set $\Hb := \{ z\in \Cb | \text{Im}\, z >0\}$, $\overline{\Hb} := \{ z\in \Cb | \text{Im}\,  z < 0\}$. 
The ground field is always assumed to be $\Cb$. 

\medskip\noindent
{\bf Acknowledgement}: I would like to thank Tom Bridgeland whose lecture delivered in Max-Planck Institute for Mathematics at Bonn in May 2007 triggered a quantum leap in my own understanding of this subject. I thank Hao Zheng for many inspiring conversations during my preparation of this paper, Fong-Ching Chen for his interesting lecture on Newton, Matilde Marcolli for clarifying my confusions on noncommutative geometry, Yi-Zhi Huang, Antun Milas, Ingo Runkel and Wang Zhong for many valuable comments on earlier drafts of this paper. I want to thank Urs Schreiber and referee whose comments on an earlier draft leads to a significant revision of it. I want to thank Arthur Greenspoon who has corrected a lot of my English grammar mistakes. I also want to thank Yi-Zhi Huang, James Lepowsky, Christoph Schweigert, Matilde Marcolli and Dennis Sullivan for their constant support for many years. This work was supported in part by the Basic Research Young Scholars Program of Tsinghua University, Tsinghua University independent research Grant No. 20101081762 and by NSFC Grant No. 20101301479.

\bigskip
\section{Geometry and physics}  \label{sec:geometry}

In this section, we will discuss the naturalness of the categorical, holographic and spectral properties of SAG from a historical, philosophical and physical point of view. The central question in this endeavor is: what is a space or geometry?

\begin{quote}
It is known that geometry assumes, as things given, both the notion of space and the first principles of constructions in space. She gives definitions of them which are merely nominal, while the true determinations appear in the form of axioms. The relation of these assumptions remains consequently
in darkness; we neither perceive whether and how far their connection is necessary, nor a priori, whether it is possible.

\hfill{--- Bernhard Riemann}
\end{quote}

The categorical nature of SAG is not only natural but also childlike. Indeed, 
our physical world as perceived from the eyes of a human being consists of various physical objects, such as trees, houses, furniture, mountains, etc, and their interrelations. Namely, our innocent view towards reality is already categorical. It is natural for us to think likewise in our search for a definition of physical space. 

Geometry is an old science grown out of the practice of ``earth-measuring". Instead of a set-theorectical foundation, geometry in ancient times regards geometrical objects, such as squares, triangles and cubes, as independent objects rather than the combinations of their points. For example, in Euclid's ``Elements", the notion of line is independent of that of point. It is even possible that the notion of point came later than that of line if they did not appear at the same time. This possibility seems natural to me because a point (or rather its approximation) rarely exists alone in nature. Instead, it appears more often as an end of a line-like object, such as a rope, or a corner of a 3-dimensional body, or an intersection of two lines.

Unfortunately, this intuitive and innocent point of view towards physical space did not lead us very far. The real revolution happened after Ren\'{e} Descartes' invention of Cartesian coordinates.  It provided a foundation of analytic geometry, in which a geometric object, such as a line or a surface, can be defined by the coordinates of its points. This is perhaps the beginning of a set-theorectial point of view of geometry, which prevailed in modern geometry in the 20th century. In particular, a space is made of points, which are ideal and abstract and have no internal structures. The content of geometry lies in how these ideal points are synthesized together. However, we cannot really glue points together due to the lack of internal structure. What one can do is only to parametrize them by known mathematical structures such as the real numbers $\Rb$ which is a mathematical model for continuum space.

The birth of calculus further strengthened this set-theorectial point of view of geometry. The success  of calculus was largely due to its applications in physics, such as in the study of the orbits of celestial bodies. The key to these applications is a deep observation that many physical processes have little to do with the internal structure of the physical objects involved in this process at least in the first order of approximation. Therefore, to study the dynamics only, one can ignore the internal structure of these physical objects. In many cases, by doing so, these physical objects can be simplified to ideal points. Their dynamics, as a consequence, can be translated into a problem of solving differential equations in calculus. This radical simplification has a deep influence on our perception of geometry. The notion of a space in geometry gradually lost its meaning as a physical space. It is treated more and more as a parameter space of a physical process or a moduli space of certain mathematical structures. For this reason, I would like to borrow terminology from computer science and call such a geometry as a {\it process oriented geometry}.

An extreme point of view along this line is to accept the notion of Cartesian spacetime as an absolute physical existence which is independent of matter living in it. This is Newton's idea of an absolute spacetime, which had a big influence on later physics and geometry. Newton's idea was not accepted without resistance. For example, Leibniz argued against it by showing that the physical existence of such absolute space is contradictory to the principle of sufficient reason and the identity of indiscernibles \cite{leibniz}. As we know, this absolute-space point of view was shattered by Einstein's theory of general relativity in which physical matter and geometry are indispensable to each other. Nevertheless, Einstein's theory, as a physical theory of geometry at large length scales, preserves many local properties of Cartesian space so that we can still apply calculus at least locally. Thus it is a global generalization of calculus. This fact is manifest in the modern definition of {\it manifold} which is locally given by open sets in $\mathbb{R}^n$. 
So it is fair to say that modern differential geometry, as a global generalization of calculus, is still a process oriented geometry. In hindsight, it is hard to tell if it is the most natural thing to do or just a historical convenience due to our strong reluctance to give up calculus as a powerful tool. Because of the success of the theory of manifolds in physics, after a few generations of physicists being educated in this framework, it is very difficult to imagine that there might be other possibilities.

Our main purpose in introducing the terminology of process oriented geometry is to introduce a different type of geometry which can very well be called {\it object oriented geometry}. In contrast to Newton, Leibniz described a space as a network of relation between physical objects, and it can not exist alone without the existence of physical objects. It is more intuitive than Newtonian space because it is exactly how we perceive the world from our daily life. But this familiar experience is not very helpful because it is too complex to be formalized. All physical objects, even atoms in modern physics, contain structures too rich to be handled easily. But complexity is not the real difficulty. There are more serious problems in this approach. Before the advent of Newtonian Gravity, a popular natural philosophy in the 17th century was the mechanical philosophy, a main contributor to which was Ren\'{e} Descartes. He believed that there is no vacuum, everything physical is made of tiny ``corpuscles" of matter, and force is possible only through the collision of adjacent corpuscles. Leibniz's network view also suggests something similar to the nature of force. From these points of view, it is hard to understand Newton's gravity because it acts at a distance. Therefore, it is a real challenge for any successful theory from Leibniz's school to understand Newton's gravity. It perhaps had to wait until modern condensed matter physics to suggest a possible answer: gravity is emergent \cite{sak}.

A real breakthrough for this object oriented point of view of geometry happened in the 1960s when Grothendieck and his school established a complete new foundation of algebraic geometry. Although the notion of scheme is still set-theoretical, Grothendieck emphasized that it should be understood as a representable functor from the category of schemes to the category of sets. In other words, a scheme can be defined by the its relation to other schemes including those that cannot be viewed as its subschemes. This is a true object oriented geometry. Even in retrospect, it is still hard to understand how Grothendieck made this breakthrough, since it was as if the historical burden had nearly no effect on him. He simply started everything from scratch. The success of Grothendieck's algebraic geometry in mathematics created a wave of replacing the set-theoretical foundation of mathematics  by a categorical one. In the 1980s and 1990s, this wave found its new driving force from another field of science: physics.

\medskip
The ultimate source of information on what a space or geometry really is should come from the observation of our physical universe. That is the field of physics. As we mentioned before, from a childlike viewpoint, the world is obviously a network of physical objects and their interactions. A physical object can only be understood through its interaction with other physical objects, beyond which there is no further reality. Only interaction is directly observable. For example, what one observes in a cloud chamber that connected to a particle accelerator is not a particle itself but its interaction with other particles. It is similar to the slogan in category theory that morphisms are more important than objects. But, for most working physicists, such a naive viewpoint is not enough to convince them to accept categories as a new universal language. They need something workable and computable. Therefore, a better way is to look at the categorical roots in modern physical theories.

Physics in the last century has gone through a quantum revolution. As we dig deeper and deeper into the microscopic world, we see a rather strange universe. Cartesian geometric intuition loses its meaning at short length scales and becomes only a low-energy illusion. It turns out that physics at the atomic scale is governed by quantum mechanics and quantum field theory (QFT). Although they are very successful theories, they are not only difficult for beginners to study but also often puzzling to experts. The key reason, as pointed out by Alain Connes, is that they lack geometric foundations. Connes went ahead to develop such a geometric foundation which is now called Noncommutative Geometry \cite{connes}. It was inspired by the Heisenberg picture of quantum mechanics. In particular, it starts with an operator algebra $\mathcal{A}$ acting on a Hilbert space $\mathcal{H}$. When enriched by a so-called Dirac operator $D$, the triple $(\mathcal{A}, \mathcal{H}, D)$ (called a spectral triple) determines a new geometry. In the case that the operator algebra $\mathcal{A}$ is commutative, the spectral triple, equipped with additional natural structures and satisfying certain natural conditions, recovers the Riemannian geometry. This is a spectral geometry by its nature.

A very different approach to quantum theory is the path integral approach. It is the most popular and productive approach to study quantum field theory. The key to its success lies in the fact that much classical geometric intuition is preserved in this formulation. Indeed, superficially, there is no non-commutative variable appearing in the path integral. Of course, it is not really classical. Its hidden noncommutativity can be seen from the fact that the free algebra of a few noncommutative variables can be viewed as coordinate functions on the path space of a lattice \cite{kapranov}. Nevertheless, the path integral approach carries a lot of commutative flavors. As we will show in Section \ref{sec:cl-cft}, a closed CFT is not commutative in the classical sense but commutative in a stringy sense. More precisely, a closed CFT is a commutative associative algebra in a braided tensor category obtained from the representation theory of a VOA. Similar commutativity also appeared in higher dimensional topological field theories \cite{lurie}. In string theory, we know that the low energy effective theory of a closed string theory reproduces Einstein's theory of gravity. I believe that this correspondence between stringy commutativity and the commutative geometry underlying Einstein's General Relativity is not accidental. Perhaps it is this hidden commutativity which forbids any naive way to quantize gravity. 

Another interesting feature of the path integral is its categorical nature. Although a proper measure for the path integral is very difficult to construct\footnote{There is a nice construction of such measure for scalar fields by Doug Pickrell \cite{doug}.}, the formal properties of the path integral are rather clear. It has been formalized by Segal \cite{Segal}, Atiyah \cite{Atiyah:1989vu} and Kontsevich as a (projective) symmetric monoidal functor from a geometric category, where the compositions of morphisms are given by gluing world-sheets, to the category of vector spaces. This gave birth to an unprecedented participation of mathematicians in the study of QFT. Category has become a common language in the study of QFT among mathematicians since then. It is important to notice that this categorical aspect of QFT is on the level of the world sheet. It turns out that the path integral also suggests to view the target space, which is our main focus, as a network of interesting subspaces connected by bordisms. This is manifest in the study of D-branes in non-linear sigma models. Geometrically, a D-brane, as a boundary condition of open strings,  can be wrapped around a subspace of the target manifold. So it can be viewed as a generalized subspace of the target manifold or simply a generalized point. What lies between two such subspaces are paths from one subspace to the other. Pictorially, we perceive the target manifold as a network of subspaces linked by paths. An algebraic model for this network is nothing but the category of D-branes introduced in Section \ref{sec:intro}. The geometric content of the space is completely encoded in this category. For example, the information of the intersection and the distance of a pair of D-branes should be extracted from the Hom space between these two D-branes. Therefore, it is reasonable to define our physical space (consisting of all observables) by the category of D-branes. The classical geometric intuition should be emergent from this mathematical structures of observables. We want to generalize this picture obtained from non-linear sigma models to all CFTs. Namely, we will forget about the target manifold and start from a closed CFT and its category of D-branes. It is just our SAG introduced in Section \ref{sec:intro}. Since the category of D-branes is itself a physical space, it is reasonable to define it to be the spectrum\footnote{One can ask why not define the spectrum to be the category of D$0$-branes. According to our philosophy of object oriented geometry, a point loses its priority as a basic building block of other objects, such as lines, surfaces and even a two-point set. Instead, a point can be viewed as an intersection of two lines for example. The relative relation between all D-branes is the full content of the geometry. Therefore, we prefer to have all D-branes in our categorical spectrum from which all kinds of moduli spaces, such as the moduli spaces of D$p$-branes for $p\geq 0$, can emerge.} of SAG. Moreover, as we discussed in Section \ref{sec:intro}, the category of D-branes is also a natural combination of the algebro-geometric spectrum (boundary conditions) and physical spectra (of Laplacians). 

Taking a closer look at this spectrum, we see that the endo-hom spaces are open CFTs compatible with a given closed CFT. We will show in Section \ref{sec:op-cl-cft} that these open CFTs are given by certain open-string vertex operator algebras. The Hom space between a pair of D-branes is given by a particular bimodule over two open-string vertex operator algebras. As we will show in later sections, open CFTs, closed CFTs and their (bi-)modules are all constructible via the representation theory of vertex operator algebra. Therefore, all these data are well-defined, constructible and computable. Moreover, sometimes these data can be reduced to more familiar mathematical structures. For example, in certain twisted theories of some supersymmetric CFTs, by considering only zero modes on D-branes,  the category of D-branes can be reduced to the bounded derived category of coherent sheaves on the target manifold \cite{sharpe}\cite{douglas}. 

The category of D-branes provides us a new way to look at geometry. From this new point of view, the classical set-theoretical geometries are not fundamental and should be viewed as emergent phenomena. For example, a set-theoretical target manifold can emerge as the moduli space of D$0$-branes \cite{aspinwall}. Moreover, restricting ourselves to zero modes on D-branes, it is known that sometimes the original target manifold (or scheme) can be reconstructed from the bounded derived category of coherent sheaves and the group of autoequivalences \cite{bondal-orlov}. This fits into our philosophy that manifolds as a concept in process oriented geometry can emerge from the moduli spaces of certain mathematical structures.

\medskip
We have seen that the categorical nature of the path integral approach to CFT naturally leads us to SAG. But perhaps there are deeper reasons for us to follow this route. In retrospect, the success of QFT is based on a rather strange logical order. Namely, we start with a very abstract notion of spacetime and rebuild everything on it. For example, the Hamiltonian is the dual of time, and fundamental particles are time-invariant states. In other words, physical matter, which lies at the heart of our sense of existence, is even a derived concept of the spacetime. Obviously, it is the Newtonian point of view lying at the core of this logic. It is rather strange from the eyes of a child or anyone belonging to the Leibniz school. Let's leave the debate between Newton and Leibniz aside. Remember, we are now looking for a theory of spacetime. It means that we have to look for some concepts that are even more basic than that of spacetime. Unfortunately, the Newtonian point of view provides no clue to such a question. But the Leibniz school of thought suggests that at least you can reverse the order of the logic by assuming that the existence of physical objects is more basic than that of spacetime. On the practical side, to work out this dual picture is certainly a way to enrich and deepen our understanding of the situation. On the philosophical side, I do believe that the question on the existence of something is deeper than the notion of spacetime. What lies at the heart of the former question is the notion of information \cite{bekenstein}, which can tell an observer that there is something. On the one hand, the total information of an observable physical object should be nicely structured. On the other hand, information itself should be some kind of structure that is observable to other structures. Therefore, we can simply define a physical object as a nicely structured stack of information. As a consequence, our physical space is nothing but a network of structured stacks of information, from which spacetime can emerge. If all the observable information is spectral, this network is already very close to our proposal of the spectrum of SAG.

\medskip
We have argued that the categorical feature of SAG is natural and in some sense a return to innocence. More importantly, it is computable. But it still might not be sufficient to convince us that it is the right track to follow. In physics, the search of a new theory is often guided by principles, which arise from experimental or theoretical observations but are excluded from existing theories. Such principles are required to lie at the core of any successful theory, and serve as severe constraints in new theories. In the search for a theory of spacetime, many physicists believe that a fundamental guiding principle is the so-called Holographic Principle\footnote{It is perhaps better to call it holographic hypothesis at the current stage.} (see \cite{bousso} for a review), which says that the boundary of the universe contains all the information of the universe.  Rooted in black hole thermodynamics, the Holographic Principle was first proposed by Gerardus 't Hooft \cite{thooft} and given a precise string-theory interpretation by Leonard Susskind \cite{susskind}. There are several different ways to formalize the idea of Holographic Principle. We will use a particular one in the context of CFT. It says that a boundary theory determines the bulk theory uniquely. Ideally, we would like to check this Holographic Principle in SAG for the most general D-branes, the so-called conformally invariant D-branes (see Definition \ref{def-ocfa-V}). Unfortunately, such a result is still not available. What we can do now is to test it in rational CFTs for D-branes respecting a chiral symmetry given by a rational VOA $V$, also called $V$-invariant D-branes (see Definition \ref{def-ocfa-V}). In these cases, the Holographic Principle indeed holds \cite{unique}\cite{cardy-alg}. More precisely, the closed CFT can be obtained by a given open CFT (over $V$) $A_\op$ by taking the center $Z(A_\op)$ of $A_\op$ (see Theorem \ref{thm:unique}). 
Moreover, such holographic phenomena also appear in string topology \cite{string-top}\cite{bct}, TCFT \cite{costello}\cite{bfn}, extended Turaev-Viro topological field theories \cite{kitaev-kong} and topological field theories in general \cite{lurie}\cite{lurie-tft}. Having partially passed the test, we conjecture that the Holographic Principle holds for many sufficiently nice conformally invariant D-branes.  In other words, our SAG should be automatically holographic. A single nice D-brane is enough to determine the bulk theory which further determines the entire geometry. In the light of William Blake's poem, we will be able to see the world in a grain of D-brane.

In physics, a strong support for the Holographic Principle comes from the AdS/CFT correspondence \cite{ads-cft} which also suggests that the Holographic Principle can provide a natural mechanism for gravity to emerge (see also for example \cite{verlinde}). This is also compatible with our philosophy. Indeed, the Holographic Principle can be viewed as only one direction of boundary-bulk duality. The inverse statement of Holographic Principle says that a bulk theory does not determine the boundary theory uniquely. The ambiguity of non-uniqueness is nothing but the spectrum of an entirely new geometry SAG, from which gravity or Riemannian geometry can emerge. 

\medskip
That SAG incorporates the Holographic Principle as an intrinsic property might not sound exciting at all to mathematicians. In mathematics, a good program is often judged by its richness and its power in solving old problems. Although one can get a glimpse of the power of CFT from its application in many topics, such as mirror symmetry, low dimensional topology and the geometric Langlands correspondence, etc, its full power in solving problems in other fields of mathematics including number theory has not been sufficiently disclosed yet. We will choose to comment on its richness instead.  

Metric is a geometric notion that plays an important role in physics. It is also considered by Riemann as a necessary ingredient of the doctrine of space. Modern algebraic geometry was developed by Grothendieck's school without any influence from physics. Its interaction with physics came rather late. This delay is partly due to the lack of metric in algebraic geometry. The recent marriage of algebraic geometry and physics via string theory perhaps is simply a sign for a possible unification of algebraic geometry with the metric. In CFT, the information of metric is encoded in the modules of the Virasoro algebra or its supersymmetric counterpart. Such modules can be roughly viewed as the spectrum of the Laplacian or Dirac operator on the loop space or path spaces. One can extract metric information from the spectrum of the Laplacian or Dirac operator similar to Connes' spectral triple. But one should keep an open mind. The intrinsic notion of metric in this 2-spectral geometry might not be the one familiar to us. One might expect a hierarchy of metrics or even more exotic notion of metric. The possible exotic behavior can be seen from our loop space intuition, which suggests that space is highly entangled at high energy through the channels of long paths. The usual intuition on distance between objects loses its pertinence at high energy. Therefore, we expect a radical new and rich notion of 2-metric in this 2-spectral geometry. In this sense, a closed CFT, as a ``stringy commutative ring" over a ``stringy ground field" constructed from the Virasoro algebra (or super-conformal algebra), provides a natural unification of algebraic geometry with the metric. 

Unification of this kind is quite inspiring if you are a believer in the unity of mathematics. In addition to it, CFT has revealed much more unifying power in mathematics. Very often, seemingly unrelated mathematical structures appear in the same CFT. For example, the first discovered CFT in mathematics was the Monster moonshine CFT \cite{FLM1}. It contains the tensor-product of 48 Ising models \cite{DMZ},  its automorphism group is the Monster group \cite{bor-va}\cite{FLM2} and its partition function is the $j$-function which is the generator of all modular functions \cite{FLM1}. What underlies all these miracles is a very deep and simple reason: QFT is a theory for physical systems with infinitely many degrees of freedom! The very existence of a mathematical structure (QFT) that is rich enough to model these infinite systems is already a miracle. The mathematical structures appearing in QFT are often infinite dimensional.  It is not surprising that one can discover many familiar but seemingly unrelated mathematical structures from a single QFT through different finite or infinite dimensional windows. The price we pay is that to construct such a rich structure of QFT explicitly is often a difficult problem.

\medskip
SAG is perhaps just what many physicists had in mind and has already been actively pursued for many years. Many aspects of SAG are known to physicists. Unfortunately, these results are scattered like leaves of a tree, but float in the air without knowing its roots. For example, the supposed-to-be-derived notions such as points, lines, paths and surfaces are unavoidably used in a lot of discussion on the emergence of spacetime in a lot of physical literature as something fundamental. How mathematicians can help is to construct the roots and provide a new calculus so that notions like points, line, surfaces, space, time and causal structure can all be derived from deeper structures and principles.

Before we leave this section, we would like to point out that  2-dimensional CFTs are certainly not the only source for such object oriented geometries. 
Notice that the key ingredients of the above SAG is the boundary-bulk duality. This duality also exists in higher dimensional TFTs perhaps in other QFTs as well. For example, for an $n$-dimensional TFT, we will have an $(n-1)$-category of boundary conditions. In many cases, boundary theories are again obtained by the internal homs of boundary conditions, and the bulk theory is obtained by taking the center of a boundary theory \cite{lurie}\cite{lurie-tft}. So I believe that it is possible to define a new geometry for many higher dimensional QFTs. The richer the structure of a QFT is, the richer the corresponding geometry is. Unfortunately, QFTs other than TFTs and 2-dimensional CFTs are still far from being understood. The structure of a CFT is certainly richer than that of a TFT. For example, there is not any spectral geometry ingredient in a TFT. Therefore, SAG based on 2-dimensional CFT is a good place to invest our labor.

\bigskip
\section{Closed conformal field theories}  \label{sec:cl-cft}

In this section, we will provide the mathematical foundations of closed CFTs without supersymmetry. From now on, by a CFT we always mean a theory without supersymmetry. A systematic study of CFT in physics start from Belavin, Polyakov and Zamolodchikov's seminal paper \cite{BPZ} in 1984. It was followed by a flood of physics papers on CFT. Some of them made deep impacts on mathematics. But the list is too long to be given here. Readers can get some physical ideas behind the materials presented in this review from \cite{BPZ}\cite{FriedanShenker}\cite{vafa}\cite{moore-seiberg}\cite{witten}. See also Francesco, Mathieu and Senechal's book \cite{fms-cft} and references therein. Now we shift our attention to the mathematical side of the story. 

\subsection{Basic definitions}

Around 1987 Segal \cite{Segal} and Kontsevich independently gave a mathematical definition of two-dimensional closed conformal field theory (closed CFT or just  CFT). This definition become very important for the development of CFT in mathematics. Segal defines a closed CFT as a projective symmetric monoidal functor $\mathcal{F}: \mathcal{RS}^b \rightarrow \mathcal{TV}$ between two symmetric monoidal categories $\mathcal{RS}^b$ and $\mathcal{TV}$. The objects of the category $\mathcal{RS}^b$ are finitely ordered sets and the morphism set $\text{Mor}(A, B)$ for $A, B\in \mathcal{RS}^b$ is the set of conformal equivalence classes of Riemann surfaces with $|A|$ ordered negatively oriented parametrized boundaries and $|B|$ ordered positively oriented parametrized boundaries, where $|A|$ and $|B|$ denote the cardinalities of the set $A$ and $B$ respectively. A negatively (positively) oriented parametrized boundary component is a boundary component equipped with a germ of conformal map from a neighborhood of the boundary to an inside (outside) neighborhood of the unit circle in complex plane such that its restriction to the boundary is real analytic\footnote{Segal's definition is more general. We only care about this restricted case here.}. The composition of morphisms is defined as sewing of Riemann surfaces along oppositely oriented parametrized boundaries. The compositions of morphisms are clearly associative. Notice that this category is not a unital category. The category $\mathcal{TV}$ is the category of locally convex complete topological vector spaces with continuous linear maps as morphisms. The projectivity of $\mathcal{F}$ means that $\mathcal{F}$, as a map on each set of morphisms, is only well-defined up to constant scalars, i.e. the image of $\mathcal{F}$ lies in the projectivization of the space of morphisms in $\mathcal{TV}$.

Segal's beautiful definition of CFT immediately leads to some interesting structures of conformal field theory such as the modular functor. However, it also has some inconveniences. 

\begin{enumerate}

\item A quantum field $\phi(x)$ in physics is defined for a point $x$ on a worldsheet. But world sheets with punctures do not live in Segal's definition. Correlation functions of CFTs in physics are usually defined as sections of certain bundles on the moduli space of Riemann surfaces with punctures \cite{FriedanShenker}. This suggests working with Riemann surfaces with punctures instead of boundaries. A geometric formulation of CFT in physics was given by Vafa \cite{vafa} as sewing of Riemann surfaces with punctures and local coordinates. 

\item Physicists usually do not work with the entire Hilbert space. Instead they often only work with a graded and dense subspace, which has already revealed a rich algebraic structure. For example, the chiral algebra (equivalent to VOA in mathematics) and its modules as crucial ingredients of a CFT are graded vector spaces. To complete properly the dense subspace to a complete topological vector space is a difficult problem. 

\end{enumerate}
In order to take advantage of the power of the algebraic method (via chiral algebras or VOAs) widely used by physicists, we would like to modify accordingly Segal's definition. 

\medskip

In physics, Vafa \cite{vafa} proposed a definition of CFT via sewing of closed Riemann surfaces with punctures after Segal. In mathematics, Igor Frenkel \cite{frenkel-ias} started a program to study CFTs in physics via the theory of vertex operator algebras even before the appearance of Segal's definition of closed CFT. One of the highlights of this program is Huang's thesis \cite{h-thesis} (see also \cite{h-gvoa}\cite{h-book}) in which Huang established the precise geometric meaning of VOA as the sewing of Riemann surfaces with punctures.   A VOA is a structure defined on a $\mathbb{Z}$-graded vector space over $\mathbb{C}$. Huang's result suggests that it is natural to consider graded vector spaces instead of complete topological vector spaces when one studies CFTs on surfaces with punctures. Vafa and Huang's works implicitly suggest a convenient and workable definition of a CFT by modifying Segal's definition only slightly. Before we do that, we need first to introduce two categories $\mathcal{RS}^p$ and $\mathcal{GV}$ below. 

\medskip
\begin{figure}
  \begin{picture}(200,200)
     \put(-50,3){\scalebox{.75}{\includegraphics{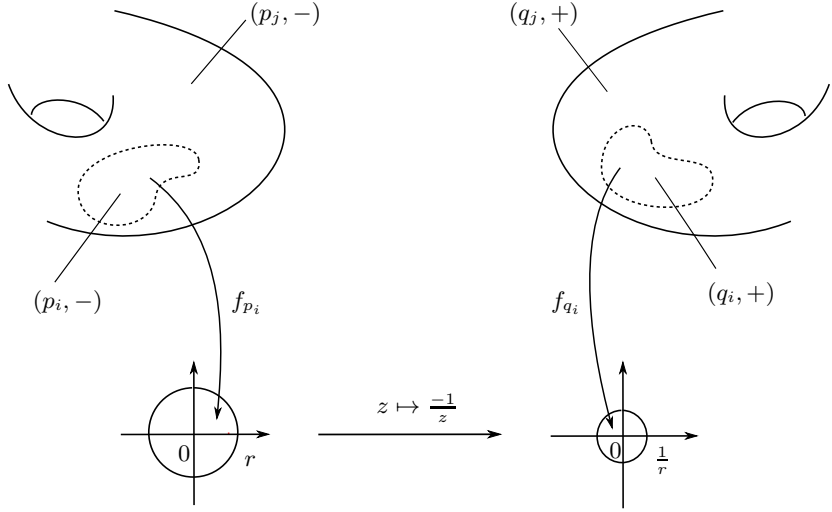}}}
   \put(90,40){$ z \mapsto \frac{-1}{z} $}
   \put(-40, 80){$ (p_i, -) $}
   \put(42, 189){$ (p_j, -) $}
   \put(215, 83){$ (q_i,+) $}
   \put(140, 189){$ (q_j, +) $}
   \put(35, 80){$ f_{p_i} $}
   \put(156, 80){$ f_{q_i} $}
   \put(40, 20){$r$}
   \put(195, 20){\small $\frac{1}{r}$} 
   \put(15, 22){$0$}
   \put(178, 23){\small $0$}
  \end{picture}
  \caption{{ \small Sewing a negatively oriented puncture $(p_i, f_{p_i}, -)$ to a positively oriented puncture $(q_i, f_{q_i}, +)$.}}
  \label{fig:sewing-operation}
\end{figure}
\noindent {\bf The category $\mathcal{RS}^p$}: 
We first introduce some terminology. By {\it a puncture in a Riemann surface} $\Sigma$, we mean a marked point $p$ in $\Sigma$ equipped with a local coordinate, which is a germ of conformal map $f_p$ from a neighborhood of $p$ to a neighborhood of $0\in \mathbb{C}$ such that $f_p(p)=0$, and an orientation, which is a label $\eps_p = \pm$ on $p$. It is called positively oriented if $\eps_p=+$, and negatively oriented if $\eps_p=-$. Such a puncture is denoted by $(p,f_p,\eps_p)$. A Riemann surface (not necessarily connected) with $k$ ordered positively oriented punctures and $l$ ordered negatively oriented punctures will be called a $(k,l)$-surface and will be denoted by
\be \label{eq:element}
(\Sigma | (p_1, f_{p_1}, +), \dots, (p_k, f_{p_k}, +); (q_1, f_{q_1},-), \dots, 
(q_l,f_{q_l}, -))~.
\ee
or $\Sigma$ for simplicity. Two $(k,l)$-surfaces are conformally equivalent if there is a biholomorphic map between them that preserves the orders, the orientations and local coordinates of the punctures. We denote the conformal equivalence class of $\Sigma$ by $[\Sigma]$ and the moduli space of $(k,l)$-surfaces of genus $g$ by $\mathcal{E}_{(k,l)}^g$. We also set $\mathcal{E}^g:= \cup_{k,l} \mathcal{E}_{(k,l)}^g$, $\mathcal{E}_{(k,l)}:=\cup_g \mathcal{E}_{(k,l)}^g$ and $\mathcal{E}:= \cup_g \mathcal{E}^g$. These moduli spaces naturally has a structure of infinity dimensional complex manifolds \cite{h-book}. 
We denote the space of germs of conformal maps $f_p$ by $F_p$. Notice that the $F_p$ are isomorphic to $F_0$ which is the space of germs of conformal maps from a neighborhood of $0\in \Cb$ to $0\in \Cb$. If we denote the moduli space of Riemann surfaces with marked points by $\mathcal{M}_{(k,l)}^g$, then $\mathcal{E}_{(k,l)}^g$ can be viewed as a bundle over $\mathcal{M}_{(k,l)}^g$ with fiber given by $F_0^{k+l}$. We set $\mathcal{M}:= \cup_{g,k,l} \mathcal{M}_{(k,l)}^g$.

Now we discuss how to sew a $(k,l)$-surface to an $(l,m)$-surface. 
Let $\Sigma_1$ be a $(k,l)$-surface with $l$ ordered negatively oriented punctures $(p_i, f_{p_i}, -), i=1, \dots, l$, and let 
$\Sigma_2$ be an $(l,m)$-surface with $l$ ordered positively oriented punctures $(q_i, f_{q_i}, +), i=1,\dots, l$. Let $D(r,0)$ denote the closed disk in $\Cb$ centered at $0$ with radius $r>0$. We say that $p_i$ can be sewn to $q_i$ for all $i=1, \dots, l$ at the same time if there exist $r_i>0, i=1,\dots, l$ such that $f_{p_i}^{(-1)}(D(r_i,0))$ in $\Sigma_1$ and $f_{q_i}^{-1}(D(1/r_i, 0))$ in $\Sigma_2$ are all disjoint. Depicted in Figure \ref{fig:sewing-operation}, the sewing operation is defined by first cutting off the interiors of the disks $f_{p_i}^{(-1)}(D(r_i,0))$ and $f_{q_i}^{-1}(D(1/r_i, 0))$, then identifying their boundaries via the maps $f_{q_i}^{-1} \circ J \circ f_{p_i}$ for $i=1,\dots, l$,  where $J: \Cb^{\times} \rightarrow \Cb^{\times}$ is given by $w \mapsto \frac{-1}{w}$. When such $r_i$ do not exist, a sewing operation cannot be defined. Notice that we defined the sewing operation by a choice of $r_i, i=1, \dots, l$. For different choices, the surfaces obtained after sewing operations are all conformally equivalent. Thus the sewing operations are well-defined on $\mathcal{E}$ (independent of choices). Notice that a sewing operation between a pair of oppositely oriented punctures can always be defined if we rescale $f_{p_i}$ or $f_{q_i}$.

Let $\mathcal{RS}^p$ be a partial category of finite ordered sets, and $\text{Mor}_{\mathcal{RS}^p}(A, B):=\mathcal{E}_{(|A|,|B|)}$ for $A, B\in \mathcal{RS}^p$, and the compositions of morphisms be given by sewing operations. The word ``partial'' refers to the fact that the sewing operations are not always well-defined. The compositions of morphisms are associative in the sense that 
$[\Sigma_1]\circ ([\Sigma_2]\circ [\Sigma_3]) = 
([\Sigma_1] \circ [\Sigma_2]) \circ [\Sigma_3]$
whenever both sides are well-defined. The symmetric monoidal structure on 
$\mathcal{RS}^p$ is given by disjoint union. It is not hard to see that $\mathcal{RS}^p$ is actually a unital category. 

\void{
\begin{rema}  {\rm 
This partial category is not satisfying. It should be possible to have a very detailed description of the meaning of ``partial-ness". Namely, we need redefine the sewing operations, and emphasis that the composition can be defined except on some singular devisors?  Perhaps analytic continuation must be used in order to define certain compositions, and one need to describe the behavior of the sewing operations near the singular points. So we expect to have a hierarchy of structures: partial CFT, smooth CFT, real analytic CFT, logarithmic CFT, meromorphic CFT, etc.  
}
\end{rema}
}

\medskip
\noindent{\bf The category $\mathcal{GV}$}:
By {\it a graded vector space} $H$ in this work, we always mean a vector space $H$ over $\Cb$ graded by an abelian group $G_H$, i.e. $H=\oplus_{n\in G_H} H_{(n)}$, with finite dimensional homogeneous spaces and a weak topology induced from the restricted dual space $H':= \oplus_{n\in G_H} H_{(n)}^*$. For the purpose of this work, the abelian group $G_H$ can be taken to be $\mathbb{Q}^i$ or $\mathbb{R}^j$ for some $i,j\in \Nb$. We will use $\overline{H}$ to denote the algebraic completion of $H$, i.e. $\overline{H}= \prod_{n\in G_H} H_{(n)}$. We define a non-associative partial category $\mathcal{GV}$. The objects in $\mathcal{GV}$ are graded vector spaces. For $H_1, H_2\in \mathcal{GV}$, the morphism set $\text{Mor}_{\mathcal{GV}}(H_1, H_2)$ is defined to be $\Hom_{\Cb}(H_1, \overline{H_2})$ instead of $\Hom_{\Cb}(H_1, H_2)$! Notice that the compositions of morphisms in $\mathcal{GV}$ are not well-defined a priori. Let $f\in \Hom_{\Cb}(H_1, \overline{H_2})$ and $g\in \Hom_{\Cb}(H_2, \overline{H_3})$. If the sum $\sum_{n\in G_{H_2}} g(P_nf(u))$, where $P_n: \overline{H_2}\rightarrow (H_2)_{(n)}$ is the projection operator, is absolutely convergent for all $u\in H_1$, then we say that $g\circ f$ is well-defined and $g\circ f(u) := \sum_{n\in G_{H_2}} g(P_nf(u))$. The composition is non-associative. For any triple of morphisms $h, g, f$, $h\circ (g \circ f) \neq (h\circ g) \circ f$ in general even if both sides are well-defined. The symmetric monoidal structure on $\mathcal{GV}$ is the usual one on the category of vector spaces.

\medskip

The moduli space $\mathcal{M}$ are naturally equipped with coordinates which can be expressed by a sequence of complex variables $(z_1,z_2,\dots)$ \cite{h-book}. Any functor defined on $\mathcal{RS}^p$ depending on each $z_i$ smoothly, holomorphically and real-analytically will be called smooth, holomorphic and real-analytic, respectively. In particular, $\mathcal{F}$ being real-analytic means that it has the following expansion:
\be  \label{eq:real-analytic}
\mathcal{F}(z_i) = \sum_{m,n\in \Rb} F_{mn} z_i^m \bar{z}_i^n, \quad\quad i=1, 2, 3, \cdots. 
\ee
For many algebraic constructions of CFTs, it is enough to assume that $\mathcal{F}$ is real-analytic. $\mathcal{F}$ is called holomorphic if $\mathcal{F}(z_i) = \sum_{n\in \Zb} F_n z_i^n$. 

\medskip
In this work, we will take the following working definition of CFT. 
\begin{defn} \label{def:cft} {\rm
A (closed) CFT is a real-analytic projective symmetric monoidal functor $\mathcal{F}: \mathcal{RS}^p \rightarrow \mathcal{GV}$. }
\end{defn}

Note that such a functor $\mathcal{F}$ requires that $\mathcal{F}([\Sigma_1]) \circ \mathcal{F}([\Sigma_2])$ is well-defined whenever the composition $[\Sigma_1] \circ [\Sigma_2]$ is well-defined in $\mathcal{RS}^p$ and 
$$
\mathcal{F}([\Sigma_1]) \circ (\mathcal{F}([\Sigma_2]) \circ \mathcal{F}([\Sigma_3])) =
(\mathcal{F}([\Sigma_1]) \circ \mathcal{F}([\Sigma_2]) ) \circ \mathcal{F}([\Sigma_3])
$$
whenever $[\Sigma_1]\circ ([\Sigma_2]\circ [\Sigma_3])$ and $([\Sigma_1] \circ [\Sigma_2]) \circ [\Sigma_3]$ are both well-defined. These are highly nontrivial conditions. For rational CFTs, one can see that the convergence domains of genus-0 correlation functions exactly match with the domains in $\mathcal{M}$ where the sewing operations are well-defined \cite{h-book}\cite{ko-ffa}. 

\begin{rema} \label{rema:real-ana}{\rm
It is possible to construct some CFTs that are smooth instead of real-analytic by some functional-analytic methods, for example via some rigorous constructions of path integral (see for example \cite{doug}). But for algebraic constructions considered in this article, it is enough to restrict to real-analytic CFTs. 
}
\end{rema}

Now we discuss a possible relation between Segal CFTs and the CFTs in Definition \ref{def:cft}. On the one hand, using parametrization, one can glue open disks containing the closed unit disk centered at $0$ or $\infty$ in $\hat{\Cb} := \Cb \cup \{ \infty \}$ to a Riemann surface with parametrized boundaries. One obtains a closed Riemann surface with ordered oriented marked points ($0$ or $\infty$ of the unit disk) and local coordinates, each of which maps a neighborhood of the point to a neighborhood of the closed unit disk centered at $0$ or $\infty$ depending on its orientation. Since the local coordinates also contain the information of orientations, we can get rid of this redundancy by requiring that all local coordinates map to neighborhoods of $0\in \Cb$. In this way, we obtain an embedding of categories $\mathcal{RS}^b \hookrightarrow \mathcal{RS}^p$. We denote the subspace of $\mathcal{E}$ consisting of elements from the morphisms in 
$\mathcal{RS}^b$ by $\mathcal{E}_{\rm Segal}$.  Therefore, from a Segal CFT one cannot automatically obtain a CFT in Definition \ref{def:cft} a priori. On the other hand, Huang proposed how to obtain a Segal CFT from a CFT in Definition \ref{def:cft} by properly completing the graded vector spaces.  Huang showed explicitly in \cite{h-funvoa1}\cite{h-funvoa2} how to do this on a substructure of genus-0 CFTs (see Theorems \ref{thm:h-voa-kalg} and \ref{thm:h-voa-segal} below).  

\medskip

\medskip
\subsection{Operads and Huang's Theorems} \label{sec:operad}

We will show how a VOA naturally appears as an ingredient of a CFT via a theorem of Huang. We will use the terminology of (partial) operad freely. The reader can find an introduction to the notion of operad in Ittay Weiss' contribution \cite{weiss} to this book.

\medskip
Let $\{ \pt \}$ be a one-point set. Since $\mathcal{F}$ is monoidal, for a finite ordered set $A\in \mathcal{RS}^p$, $\mathcal{F}(A)\cong \mathcal{F}(\{ \pt \})^{\otimes |A|}$. For a morphism $f\in \text{Mor}_{\mathcal{RS}^p}(A, B)$, 
$\mathcal{F}(f)$ can be viewed as a morphism $\mathcal{F}(\{ \pt \})^{\otimes |A|} \rightarrow \overline{\mathcal{F}(\{ \pt \})^{\otimes |B|}}$. Therefore, the structure of a CFT can be transported to a structure on $\mathcal{F}(\{ \pt \})$. Namely, a CFT is nothing but a graded vector space $$V_\cl:=\mathcal{F}(\{ \pt \})$$ equipped with multilinear maps $\mathcal{F}([\Sigma]): V_\cl^{\otimes k} \rightarrow \overline{V_\cl^{\otimes l}}$ for all $(k,l)$-surfaces $\Sigma$ and $k,l\in \Nb$ satisfying all the conditions required by the functoriality of $\Fc$. For this reason, we will also denote a CFT by 
$$
(V_\cl, \{ \mathcal{F}([\Sigma])\}_{[\Sigma]\in \mathcal{E}}).
$$ 
This notation $(V_\cl, \{ \mathcal{F}([\Sigma])\}_{[\Sigma]\in \mathcal{E}})$ for CFT is very convenient for us to address substructures of a CFT. For example, the space $\mathcal{E}$ equipped with sewing operations is just a partial PROP. The CFT structure $(H, \{ \mathcal{F}([\Sigma])\}_{[\Sigma]\in \mathcal{E}})$ is nothing but a projective algebra over this partial PROP.  The sewing operations restricted to $\mathcal{E}^0$ induce the structure of a partial dioperad \cite{gan} on $\mathcal{E}^0$ called the sphere partial dioperad in \cite{cardy-cond}. Therefore, a genus-0 CFT, or equivalently the structure $(V_\cl, \{ \mathcal{F}([\Sigma])\}_{[\Sigma]\in \mathcal{E}^0})$ where all sewing operations resulting in higher genus surfaces are not allowed, is nothing but a smooth projective algebra over the sphere partial dioperad $\mathcal{E}^0$. A genus-0,1 CFT is defined to be the structure $(V_\cl, \{ \mathcal{F}([\Sigma])\}_{\Sigma\in \mathcal{E}^0\cup \mathcal{E}^1})$. We will also be interested in the space $K:=\cup_{k\in\Nb} \mathcal{E}_{(k, 1)}^0$ which has the structure of a partial operad called the sphere partial operad in \cite{h-book}. Then the structure $(V_\cl, \{ \mathcal{F}([\Sigma])\}_{[\Sigma]\in K})$ is nothing but a smooth projective $K$-algebra. 

\medskip
With all this terminology, we can state Huang's fundamental result. 

\begin{thm}[\cite{h-book}]\label{thm:h-voa-kalg}
A VOA $V$ canonically gives a holomorphic projective $K$-algebra
or equivalently a structure $(V_\cl, \{ \mathcal{F}([\Sigma])\}_{[\Sigma]\in K})$ for $V_\cl=V$ as a substructure of a CFT. 
\end{thm}

 We will illustrate the main structures of a VOA in Section \ref{sec:voa} by the above theorem. Namely, a VOA is a vector space $V$ equipped with multilinear maps $\{ \Fc([\Sigma]): V^n \to \overline{V} | \Sigma \in \Ec_{(n,1)}^0 \}_{n=0}^{\infty}$ satisfying conditions required by the functoriality of $\Fc$. 

\begin{rema} {\rm
Huang's original theorem says that the category of VOAs with central charge $c$ is isomorphic to the category of holomorphic algebras (not projective!) over a certain partial-operad extension of $K$ involving the determinant line bundle over $K$ and satisfying additional natural properties. Therefore, VOA should be viewed as a natural ingredient of holomorphic CFT with natural properties. To avoid technical details, we are satisfied with the weak version stated in Theorem \ref{thm:h-voa-kalg} for the purpose of this paper. For more details, one should consult Huang's book \cite{h-book}. 
}
\end{rema}

Similarly, we can also denote a Segal CFT by $(H_{\rm Segal}, \{ \mathcal{F}_{\rm Segal} ([\Sigma]) \}_{[\Sigma]\in \mathcal{E}_{\rm Segal}})$, a genus-0 Segal CFT as 
$(H_{\rm Segal}, \{ \mathcal{F}_{\rm Segal} ([\Sigma]) \}_{[\Sigma]\in \mathcal{E}_{\rm Segal}^0})$ where $H_{\rm Segal}$ is a locally convex complete topological vector space and $\mathcal{E}_{\rm Segal}^0 := \mathcal{E}_{\rm Segal} \cap \mathcal{E}^0$. Let $K_{\rm Segal} := \mathcal{E}_{\rm Segal} \cap K$. Sewing operations are 
always defined on $K_{\rm Segal}$. Thus $K_{\rm Segal}$ is an operad instead of a partial operad. The structure $(H_{\rm Segal}, \{ \mathcal{F}_{\rm Segal} ([\Sigma]) \}_{[\Sigma]\in K_{\rm Segal}})$ is nothing but  a smooth $K_{\rm Segal}$-aglebra structure on $H_{\rm Segal}$.  Huang also showed how to obtain a smooth $K_{\rm Segal}$-aglebra by properly completing a finitely generated VOA.
\begin{thm}[\cite{h-funvoa1}\cite{h-funvoa2}] \label{thm:h-voa-segal}
By properly completing a finitely generated VOA $V$ to a locally convex complete topological vector space $\tilde{V}$, one can obtain a holomorphic projective $K_{Segal}$-algebra on $H_{\rm Segal} = \tilde{V}$. 
\end{thm}

\begin{rema} {\rm 
Note that $V\subset \tilde{V} \subset \overline{V}$. The topology on $\tilde{V}$ in Huang's construction is even nuclear \cite{cosc}. Such a  projective $K_{Segal}$-algebra $(H_{\rm Segal}, \{ \mathcal{F}_{\rm Segal} ([\Sigma]) \}_{[\Sigma]\in K_{\rm Segal}})$ gives a substructure of a Segal CFT. 
}
\end{rema}

\subsection{Vertex operator algebras}  \label{sec:voa}

In this paper, we will not give the standard definition of VOA in terms of formal variables because it is not very illuminating for our purpose. Intrigued readers can consult many reviews of VOAs (see for example \cite{lepowsky-li}). Instead, we will introduce all ingredients of a VOA $V_\cl=V$ from the structure $(V, \{ \mathcal{F}([\Sigma])\}_{\Sigma\in K})$ of a holomorphic projective $K$-algebra by Huang's theorem \ref{thm:h-voa-kalg}, and summarize it into a working definition of VOA (see Definition \ref{def:voa}).
To avoid the technicalities associated to the word ``projective'', let us pretend that a VOA $V$ gives a $K$-algebra instead of a projective $K$-algebra. 

\medskip
Building blocks of a VOA: 

\begin{enumerate}

\item {\it Grading properties}: A VOA $V$ is a $\mathbb{Z}$-graded vector space such that 
$\dim V_{(n)} < \infty$ and the grading is bounded from below, i.e. 
$V=\oplus_{n\in \Zb} V_{(n)}$ and $V_{(n)}=0$ for $n<<0$. For $u\in V_{(n)}$, $n$ is also called {\it conformal weight} of $u$, denoted by $\text{wt} \, u$. The restricted dual space is defined by $V':=\oplus_n (V_{(n)})^\ast$. 

\item {\it Unit}: Recall (\ref{eq:element}) and consider the element $\Sigma_\one:=[(\hat{\Cb}, (\infty, f_{\infty},-))]\in \mathcal{E}_{(0,1)}^0\subset K$ where $f_{\infty}: w\mapsto \frac{-1}{w}$. This element\footnote{We choose to use a convention that is slightly different from the one used in Huang's book \cite{h-book}, where $f_\infty$ is chosen to be $w\mapsto \frac{1}{w}$.} in $K$ can be sewn with any positively oriented puncture on any sphere in $K$ and the surface thus obtained is simply a sphere with one fewer positively oriented punctures but with all other data unchanged. This gives rise to a distinguished element $\mathcal{F}([(\hat{\Cb}, (\infty, f_{\infty},-))])$ in $\overline{V}$, denoted by $\one$. 
In the case of VOA, $\one \in V_{(0)}$. 

\item In the case of VOA, the $\mathcal{F}$ image of the element $\Sigma_\id:= [(\hat{\Cb},(0, f_0, +), (\infty, f_{\infty},-))] \in \mathcal{E}_{(1,1)}^0$, where $f_0: w\rightarrow w$ and $f_{\infty}: w\mapsto \frac{-1}{w}$, is nothing but $\id_V$. Moreover, for $f_0: w\rightarrow a w, \, a\in \Cb^\times$ and $f_{\infty}: w\mapsto \frac{-1}{w}$, $\mathcal{F}[(\hat{\Cb},(0, f_0, +), (\infty, f_{\infty},-))]$ acts on $V_{(n)}$ as $a^{-n}$.  

\item {\it Vertex operator}: Consider the following element: 
$$
P(z):= [(\hat{\Cb}, (z, f_z, +), (0, f_0, +), (\infty, f_{\infty},-))] \in 
\mathcal{E}_{(2,1)}^0
$$
where $z\in \Cb^{\times}, f_z: w\rightarrow w-z, f_0: w\rightarrow w, f_{\infty}: w\rightarrow \frac{-1}{w}$.  The linear map $\mathcal{F}(P(z))$ gives the vertex operator: 
$Y(\cdot, z)\cdot := \mathcal{F}(P(z)): V \otimes V \rightarrow \overline{V}$ 
for $z\in \Cb^{\times}$. For $u,v\in V$, we have 
$$
u\otimes v \mapsto Y(u,z)v = \sum_{n\in \Zb} u_n v z^{-n-1},
$$
where $u_n: V \rightarrow \overline{V}$. In the case of VOA, $u_n\in \text{End}(V)$. Moreover, if $u\in V_{(m)}$, then 
$u_n: V_{(k)}\rightarrow V_{(k+m-n-1)}$. Therefore, $u_n v =0$ for $n>>0$. 

\item {\it Unit properties}: If we sew the puncture at $\infty$ on $[\Sigma_\one]$ to the puncture at $z$ on $P(z)$, we obtain the sphere $[\Sigma_\id]$. This sewing identity is transported by $\Fc$ to the identity $Y(\one, z) =\id_V$, which is called the {\it left unit property} of the VOA. If we first sew the puncture on $[\Sigma_\one]$ to the puncture at $0$ in $P(z)$, then take the limit $z \to 0$, we obtain $[\Sigma_\id]$ again. Using $\Fc$, it gives the identity: $\lim_{z\to 0} Y(\cdot, z) \one =  \id_V$, which is called the {\it right unit property} of the VOA. 

\item {\it Convergence property I}: The puncture at $0$ on $P(z_1)$ can be sewn to that at $\infty$ on $P(z_2)$ if and only if $|z_1|>|z_2|>0$. By the axioms of CFT, this sewability implies the following convergence property: the sum
$$
\langle u', Y(u_1, z_1)Y(u_2, z_2)u_3\rangle := 
\sum_{n\in \Zb} \langle u', Y(u_1, z_1)P_nY(u_2, z_2)u_3\rangle
$$ 
is absolutely convergent when $|z_1|>|z_2|>0$ for all $u_1,u_2,u_3\in V$ and $u'\in V'$. Actually, a VOA satisfies a stronger {\it convergence property I}: the sum
\bea
&&\langle u', Y(u_1, z_1)\cdots Y(u_n, z_n)v\rangle   \nn
&&\hspace{2cm} =
\sum_{k_1,\cdots, k_n \in \Zb} \langle u', Y(u_1, z_1)P_{k_1}Y(u_2, z_2) \cdots P_{k_n} Y(u_n, z_n)v\rangle
\nonumber
\eea
is absolutely convergent when $|z_1|> \cdots >|z_n|>0$ for all $u_1,\cdots, u_n, v\in V$ and $u'\in V'$.

\item {\it Convergence property II}: The puncture at $z_2$ on $P(z_2)$ can be sewn to 
that at $\infty$ on $P(z_1-z_2)$ if and only if $|z_2|>|z_1-z_2|>0$. This sewability
implies the following convergence property: the sum
$$
\hspace{1cm} \langle u', Y(Y(u_1, z_1-z_2)u_2, z_2) u_3\rangle := 
\sum_{n\in \Zb} \langle u', Y(P_nY(u_1, z_1-z_2)u_2, z_2) u_3\rangle
$$
is absolutely convergent when $|z_2|>|z_1-z_2|>0$ all $u_1,u_2,u_3\in V$ and 
$u'\in V'$. 

\item {\it Associativity}: The two sewing operations defined in the convergence property I and II are both well-defined when $|z_1|>|z_2|>|z_1-z_2|>0$. They give the same surface. By applying $\mathcal{F}$, we obtain the associativity of VOA: 
$$
\langle u', Y(u_1, z_1)Y(u_2, z_2)u_3\rangle =
\langle u', Y(Y(u_1, z_1-z_2)u_2, z_2) u_3\rangle 
$$
when $|z_1|>|z_2|>|z_1-z_2|>0$ for all $u_1,u_2,u_3\in V$ and $u'\in V'$. It is nothing but the operator product expansion (OPE) in physics.

\item {\it Commutativity}: By sewing the $0$ on $P(z_1)$ to $\infty$ on $P(z_2)$ when $|z_1|>|z_2|>0$ and sewing the $0$ on $P(z_2)$ to $\infty$ on $P(z_1)$ when $|z_2|>|z_1|>0$, we obtain two different elements in $\mathcal{E}_{(3,1)}^0$ which are path connected. 
By the fact that $\mathcal{F}$ is holomorphic and monodromy-free\footnote{$\mathcal{F}$ is well-defined on the space $\mathcal{E}_{(3,1)}^0$, which is not simply connected.}, we see that the $\mathcal{F}$ images of these two elements in 
$\mathcal{E}_{(3,1)}^0$ must be the unique analytic continuations of each other. 
This implies that 
$$
\langle u', Y(u_1, z_1)Y(u_2, z_2)u_3\rangle, \quad\quad\quad |z_1|>|z_2|>0,
$$
and 
$$
\langle u', Y(u_2, z_2)Y(u_1, z_1)u_3\rangle,  \quad\quad\quad |z_2|>|z_1|>0,
$$
are analytic continuations of each other for all $u_1,u_2,u_3\in V$ and $u'\in V'$. 
This property is called the commutativity or the locality of VOA. 

\item {\it Conformal properties}: Consider the element $\Sigma_{\omega}:=[(\hat{\Cb}, (\infty, f_{\infty}, -))]\in \mathcal{M}_{(0,1)}^0$ with $f_{\infty}: w\rightarrow \frac{-1}{w-\eps}$. It is clear that $\mathcal{F}(\Sigma_{\omega})$ is an element in $\overline{V}$. Then 
$\omega:=\frac{d}{d\eps}|_{\eps=0} \mathcal{F}(\Sigma_{\omega})$ gives a distinguished element in $\overline{V}$. In the case of VOA, $\omega\in V_{(2)}$ and we have 
\be  \label{eq:omega}
Y(\omega, z) = \sum_{n\in \Zb} L(n) z^{-n-2}~,
\ee
where $L(n), n\in \Zb$, can be shown by the sewing properties of spheres in $\Ec_{(1,1)}$, to generate a Virasoro Lie algebra (see the proof of \cite[Prop. 5.4.4]{h-book}), i.e. 
\be  \label{eq:Virasoro}
[L(m), L(n)] = (m-n)L(m+n) + \frac{c}{12}(m^3-m) \delta_{m+n, 0}~, \quad\quad \forall m,n\in \Zb,
\ee
where $c\in \Cb$ is called the central charge\footnote{We cheat here. This is the place we cannot ignore the ``projectivity". In order to obtain the Virasoro algebra with a nontrivial central charge, we need extend the sphere partial operad $K$ by the $\frac{c}{2}$-power of determinant line bundle over $K$ (see \cite{h-book} for more details).}. Moreover, $\omega=L(-2)\one$ and $L(0)$ acts on $V$ as the grading operator. We also have 
\be   \label{eq:L(-1)}
[L(-1), Y(u, z)] = Y(L(-1)u, z) = \frac{d}{dz} Y(u, z)
\ee 
which can be easily derived by considering sewing the puncture on $\Sigma_\omega$ to the puncture at $z$ on $P(z)$. Using (\ref{eq:L(-1)}) or sewing the puncture on $\Sigma_\one$ to the puncture at $0$ on $P(z)$, one can derived the following identity: 
$$
Y(u, z) \one = e^{zL(-1)} u
$$
which is sometimes called the creation property. 

\end{enumerate}

We can summarize some of the properties of VOA listed above into a working definition of VOA: Definition \ref{def:voa}. This definition is equivalent to the usual definition of VOA \cite[Prop. 1.7]{h-jacob}.
\begin{defn}   \label{def:voa} {\rm
A VOA is a quadruple $(V, Y, \one, \omega)$ where $V=\oplus_{n\in \Zb}V_{(n)}$ is a $\Zb$-graded vector space with grading operator $L(0)$, $\one \in V_{(0)}$, $\omega\in V_{(2)}$ and $Y$ is a vertex operator: 
\bea
Y(\cdot, z)\cdot: V\otimes V &\to& \overline{V}   \nn
u \otimes v &\mapsto& Y(u,z)v = \sum_n u_n v z^{-n-1},  \quad u_n \in \text{End} V, \nonumber
\eea
satisfying the grading properties, the unit properties, the convergence property I and II, the associativity and the commutativity listed above and the identity (\ref{eq:omega}), (\ref{eq:Virasoro}), (\ref{eq:L(-1)}). 
}
\end{defn}

It is straight forward to define the notion of a module over a VOA. We will not do it here. Readers can consult \cite{lepowsky-li}. The tensor product theory of modules over a VOA is rather complicated. A relative easier notion is the so-called intertwining operator \cite{FHL}, which is used to define the tensor product. We will not recall this notion here either. Readers who are interested in these topics should consult \cite{FHL}\cite{hl-tcat-1}\cite{hl-tcat-2}\cite{hl-tcat-3}\cite{h-tcat-4}.

Notice that we have only used spheres $\Sigma_\one, \Sigma_\id,  P(z), \Sigma_\omega$ and spheres in $\Ec_{(1,1)}$ to define the data in a VOA. But these spheres generate the entire sphere partial operad $K$ by sewing operations. Therefore, the data in a VOA is enough to construct a projective $K$-algebra $(V, \{ \mathcal{F}([\Sigma])\}_{\Sigma\in K})$. An explicit formula for $\Fc([\Sigma])$ for a generic element $[\Sigma] \in K$ in terms of the data of a VOA is known \cite[eq. (5.4.1)]{h-book}. Intrigued readers should consult Huang's book \cite{h-book}.

\subsection{Holomorphic CFTs}

Since one cannot obtain all elements in $\mathcal{E}^0$ from those in $K$ by sewing operations,  a VOA is not enough to construct a genus-0 CFT.  It was proved in \cite{ko-ffa} that a VOA $V$ equipped with a non-degenerate invariant bilinear form on $V$ is enough to give a holomorphic genus-0 closed CFT $(V_\cl, \{ \mathcal{F}([\Sigma])\}_{\Sigma\in \mathcal{E}^0})$ by taking $V_\cl=V$. The notion of an invariant bilinear form on a VOA was first introduced in \cite{FHL}. More precisely, a bilinear form $(\cdot, \cdot): V\otimes V \to \Cb$ is called invariant if the following identity:
$$
(w, Y(u, z)v) = (Y(e^{-zL(1)} z^{-2L(0)}u, -z^{-1})w, v)
$$
holds for $u,v,w\in V$ and $z\in \Cb^\times$. This definition of invariant bilinear form is slightly different from that given in \cite{FHL}. We follow the convention taken in \cite{ko-ffa} in order to obtain simpler categorical formulation later. Such an invariant bilinear form is automatically symmetric. For simplicity, a VOA equipped with a non-degenerate invariant bilinear form will be called a self-dual VOA. 

Note that all elements in $\mathcal{E}$ can be obtained by applying sewing operations on elements in $\mathcal{E}^0$ repeatedly. Therefore, if a genus-0 CFT is extendable to all higher genera, then it uniquely determines the entire CFT. Higher genus theories provide no new data but only some compatibility conditions. Therefore, Huang's theorem \ref{thm:h-voa-kalg} suggests that one might be able to obtain a CFT from a self-dual VOA. If this is indeed possible, then this CFT will be holomorphic. A CFT requires $\mathcal{F}$ to be well-defined on $\mathcal{E}$. Equivalently, $\mathcal{F}$ is invariant under the actions of all mapping class groups. In the case of holomorphic CFTs, this invariance property puts serious constrain on the VOAs. So far, it is not known which VOAs satisfy all the necessary properties. But if we restrict to genus-0,1 CFT, the answer is known. The examples of such VOAs are the Monstrous Moonshine VOA $V^{\natural}$ \cite{FLM1}\cite{bor-va}\cite{FLM2} and the VOAs constructed from self-dual positive even lattices \cite{FLM2}. Since such VOAs are candidates to construct holomorphic CFTs, they are called holomorphic VOAs. 

\subsection{Full field algebras}  \label{sec:ffa}

In general, we cannot expect CFTs to be holomorphic. In this work, we are interested in real-analytic CFTs (recall Remark \ref{rema:real-ana}). Again, we first restrict our attention to $K$. Can we obtain a real-analytic $K$-algebra from VOAs? The answer is yes. Given a pair of VOAs $(V^L, Y_{V^L}, \one_{V^L},\omega_{V^L})$ and $(V^R, Y_{V^R}, \one_{V^R}, \omega_{V^R})$, one can combine them in a way to yield another interesting structure. More precisely, we set $V_\cl=V^L\otimes V^R$, 
$$
\one_\cl = \one_{V^L} \otimes \one_{V^R}, \quad\quad \omega_\cl = \omega_{V^L} \otimes 
\one_{V^R} + \one_{V^L}\otimes \omega_{V^R}, 
$$
and 
\be  \label{eq:ffa-Y}
Y_\cl( u^L\otimes u^R; z, \bar{z}) (v^L\otimes v^R) = 
Y_{V^L}(u^L, z)v^L \otimes Y_{V^R}(u^R, \bar{z})v^R~,
\ee
for $u^L, v^L\in V^L, u^R, v^R\in V^R$. Such a quadruple $(V_\cl, Y_\cl, \one_\cl, \omega_\cl)$ is not a VOA, but a so-called full field algebra introduced in \cite{ffa}. This quadruple is not much different from a VOA. In particular, if one replaces the $\bar{z}$ in (\ref{eq:ffa-Y}) by $z$, one obtains again a VOA \cite{FHL}. Thus the tensor product $V^L\otimes V^R$ has both a VOA structure and a full field algebra structure depending on how the vertex operator is defined. The latter structure is interesting because the quadruple $(V_\cl, Y_\cl, \one_\cl, \omega_\cl)$ gives another projective $K$-algebra which is not holomorphic but real-analytic. A full field algebra can be viewed as a real-analytic analogue of a VOA.

In general, if $V^L$ and $V^R$ are not holomorphic VOAs, one cannot construct the entire CFT structure on a full field algebra $V^L\otimes V^R$ due to the lack of the modular invariance property on $V^L\otimes V^R$. However, if a VOA $V$ is rational\footnote{$V$ is rational if it satisfies the conditions in Theorem \ref{thm:huang-mtc}.}, it was proved by Huang (see Theorem \ref{thm:huang-mtc}) that $\mathcal{C}_V$, the category of $V$-modules, is a modular tensor category \cite{retu}\cite{turaev-bk} on which all the mapping class groups act. If $V^L$ and $V^R$ are rational, then the VOA $V^L\otimes V^R$ is also rational \cite{DMZ}\cite{ffa}. In this case, it is reasonable to expect that we might be able to obtain a CFT of all genera by extending the full field algebra $V^L\otimes V^R$ by adding modules over $V^L\otimes V^R$ viewed as a VOA. Namely, we will look for a CFT extension $V_\cl$ of $V^L\otimes V^R$ as a $V^L\otimes V^R$-module.  Such a CFT, if it exists, will be called a CFT over $V^L\otimes V^R$. We will denote the canonical embedding $V^L \otimes V^R \hookrightarrow V_\cl$ by $\iota_\cl$. In this case,  we can use the powerful tools of tensor category to give a classification of CFTs over $V^L\otimes V^R$. A full field algebra extension of $V^L\otimes V^R$ is called a full field algebra over $V^L\otimes V^R$.

We will not give a precise definition of a full field algebra. Being a projective $K$-algebra, a full field algebra has properties similar to those of a VOA. In particular, it also satisfies a certain associativity and a commutativity. The difference lies only in their analyticity properties. We will describe only a few crucial ingredients and properties of a full field algebra that are important for the purpose of this paper. See \cite{ffa}\cite{ko-ffa} for more details. 


\bnu

\item {\it Grading properties}: $V_\cl$ is a $\Rb \times \Rb$-graded vector space, i.e. $V_\cl = \oplus_{m,n \in \Rb} (V_\cl)_{(m,n)}$, such that $\dim V_{(m,n)} < \infty$ and $V_{(m,n)}=0$ if $m-n \notin \Zb$ or $m<<0$ or $n<<0$.

\item {\it Unit}: $\one_\cl \in (V_\cl)_{(0, 0)}$. 

\item {\it Vertex operator}: for $u, v\in V_\cl$ and $z, \zeta \in \Cb^\times$ with branching cut: $-\pi <\text{Arg}(z) \leq \pi$, $-\pi \leq \text{Arg}(\zeta) < \pi$; 
\bea
Y_\cl(\cdot; z, \zeta) : V_\cl \otimes V_\cl &\to& \overline{V_\cl}  \nn
u \otimes v &\mapsto& Y_\cl(u; z, \zeta)v = \sum_{m,n\in \Rb} u_{(m,n)}v z^{-m-1} \zeta^{-n-1},   \nonumber
\eea
where for $u\in (V_\cl)_{(a,b)}$, $u_{(m,n)}: (V_\cl)_{(k,l)} \to (V_\cl)_{(k+a-m-1,l+b-n-1)}$.


\item {\it Unit property}: $Y_\cl(\one; z, \zeta) = \id_{V_\cl}$. 

\item {\it Convergence property I}: Let $V_\cl'$ be the restricted dual space of $V_\cl$, i.e. $V_\cl' = \oplus_{m,n \in \Rb} (V_\cl)_{(m,n)}^\ast$. For $u_1, u_2, v\in V_\cl, v'\in V_\cl'$, the following sum
$$
\hspace{1.3cm}\langle v', Y_\cl(u_1; z_1, \zeta_1) Y_\cl(u_2; z_2, \zeta_2)v\rangle  
= \sum_{m,n} \langle v', Y_\cl(u_1; z_1, \zeta_1) P_{(m,n)} Y_\cl(u_2; z_2, \zeta_2)v\rangle
$$
where $P_{(m,n)}: \overline{V_\cl} \to  (V_\cl)_{(m,n)}$ is the projection operator, is absolutely convergent when $|z_1|>|z_2|>0$ and $|\zeta_1|>|\zeta_2|>0$.

\item {\it Convergence property II}: For $u_1, u_2, v\in V_\cl, v'\in V_\cl'$, the following sum
\bea
&&\langle v', Y_\cl(Y_\cl(u_1; z_1-z_2, \zeta_1-\zeta_2)u_2; z_2, \zeta_2)\rangle   \nn
&&\hspace{2cm}= \sum_{m,n} \langle v', Y_\cl(P_{(m,n)}Y_\cl(u_1; z_1-z_2, \zeta_1-\zeta_2)u_2; z_2, \zeta_2)\rangle
\nonumber
\eea
is absolutely convergent when $|z_2|>|z_1-z_2|>0$ and $|\zeta_2|>|\zeta_1-\zeta_2|>0$.

\item {\it Associativity}: For $u_1, u_2, v\in V_\cl, v'\in V_\cl'$, 
\bea
&&\langle v', Y_\cl(u_1; z_1, \bar{z}_1) Y_\cl(u_2; z_2, \bar{z}_2)v\rangle  \nn
&&\hspace{2cm}=\langle v', Y_\cl(Y_\cl(u_1; z_1-z_2, \bar{z}_1-\bar{z}_2)u_2; z_2, \bar{z}_2)\rangle 
\nonumber
\eea
when $|z_1|>|z_2|>|z_1-z_2|>0$. This associativity is also called the OPE of bulk fields in physics.

\item {\it Commutativity}: For $u_1, u_2, v\in V_\cl, v'\in V_\cl'$
$$
\langle v', Y_\cl(u_1; z_1, \zeta_1) Y_\cl(u_2; z_2, \zeta_2)v\rangle  
$$
in the domain $\{ |z_1|>|z_2|>0, |\zeta_1|>|\zeta_2|>0 \}$ is the analytic continuation of
$$
\langle v',  Y_\cl(u_2; z_2, \zeta_2)Y_\cl(u_1; z_1, \zeta_1)v\rangle
$$
in the domain $\{ |z_2|>|z_1|>0, |\zeta_2|>|\zeta_1|>0 \}$ along a path of $z_2$ (fixing $z_1$) running clockwisely around $z_1$ with zero winding number and a path of $\zeta_2$ (fixing $\zeta_1$) running counter-clockwisely around $\zeta_1$ with zero winding number. Both paths do not cross the branching cuts.

\item {\it Left Virasoro element}: $\omega_L \in (V_\cl)_{(2,0)}$; {\it right Virasoro element}: $\omega_R \in (V_\cl)_{(0,2)}$. We have 
\bea
Y_\cl(\omega_L; z, \zeta) &=& \sum_{n\in \Zb} L^L(n) z^{-n-1},   \nn
Y_\cl(\omega_R; z, \zeta) &=& \sum_{n\in \Zb} L^R(n) \zeta^{-n-1},  \nonumber
\eea
where $\{ L^L(n) \}$ generate a Virasoro algebra of central charge $c^L \in \Cb$ and 
$\{ L^R(n) \}$ generate a Virasoro algebra of central charge $c^R\in \Cb$. $[L^L(m), L^R(n)]=0$ for $m,n\in \Zb$. Moreover, $L^L(0)$ and $L^R(0)$ are left and right grading operators, i.e. $L^L(0)|_{(V_\cl)_{(m,n)}} = m\, \id_{(V_\cl)_{(m,n)}}$ and $L^R(0)|_{(V_\cl)_{(m,n)}} = n\, \id_{(V_\cl)_{(m,n)}}$. 

\enu

\begin{rema}  \label{rema:comm-ffa} {\rm 
Among of all these properties, the commutativity is especially interesting to us. It is nothing but the stringy commutativity of a closed CFT mentioned in the Introduction. In Section \ref{sec:cat-ffa}, we will give a categorial formulation of a full field algebra over $V^L\otimes V^R$ where both $V^L$ and $V^R$ are rational. In this categorical formulation, the commutativity listed above reduces to the usual commutativity of an algebra in a braided tensor category \cite{ko-ffa}(see Theorem \ref{thm:genus-0-cft}). 
}
\end{rema}

For readers who are interested in how a full field algebra can produce a projective $K$-algebra, a construction and a detailed proof can be found in \cite[Thm. 1.19]{ko-ffa}. To lift a projective K-algebra $V_\cl$ to a genus-0 CFT, we need a bilinear form on $V_\cl$. Let's define $Y_\cl^f(u; x, y)v = \sum_{(m,n)} u_{(m,n)}v x^{-n-1}y^{-n-1}$ where $x, y$ are formal variables. A bilinear form $(\cdot, \cdot)_\cl: V_\cl\otimes V_\cl \to \Cb$ is called invariant \cite[eq. (2.1)]{ko-ffa}) if the following identity: 
$$
(w, Y_\cl^f(u; x, y)v)_\cl =
(Y_\cl^f(e^{-xL^L(1)-yL^R(1)} x^{-2L^L(0)}y^{-2L^R(0)}u, e^{\pi i} x^{-1}, e^{-\pi i} y^{-1})w, v)_\cl 
$$
holds for $u, v, w\in V_\cl$. Such an invariant bilinear form is automatically symmetric \cite[Prop. 2.3]{ko-ffa}. A full field algebra equipped with a non-degenerate invariant bilinear form is called self-dual. It was proved in \cite[Thm. 2.7]{ko-ffa} that a self-dual full field algebra canonically gives a genus-0 CFT.

\medskip
\subsection{Modular tensor categories}

From now on, we will assume that both VOAs $V^L$ and $V^R$ are rational. We will first recall the notion of a modular tensor category \cite{retu}\cite{turaev-bk}. 

Let $\mathcal{C}$ be a tensor category. We set the notations: $\one$ for the unit object; $\otimes: \mathcal{C} \times \mathcal{C}\rightarrow \mathcal{C}$ for the tensor bifunctor; $l_A: \one \otimes A \cong A$ and $r_A: A\otimes \one \cong A$ for $A\in \mathcal{C}$ for unit isomorphisms; $\alpha_{A,B,C}: (A\otimes B) \otimes C \cong A\otimes (B\otimes C), \forall A,B,C\in \mathcal{C}$ for the associator. If $\mathcal{C}$ is braided, we use $c_{A,B}: A\otimes B \rightarrow B\otimes A$, 
for $A,B\in \Cc$, to denote the braiding isomorphisms. 

If $\mathcal{C}$ is rigid, each object $U$ is equipped with a left dual $ ^{\vee}U$ and a right dual $U^{\vee}$. A ribbon category is a rigid braided tensor category with a twist isomorphism $\theta_U: U \rightarrow U$ for each $U\in \mathcal{C}$ satisfying certain balancing properties. In particular, in a ribbon category, one can take $^{\vee}U=U^{\vee}$. In this case, we will write the dualities as
$$
\begin{array}{llll}
  \raisebox{-8pt}{
  \begin{picture}(26,22)
   \put(0,6){\scalebox{.75}{\includegraphics{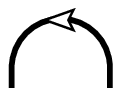}}}
   \put(0,6){
     \setlength{\unitlength}{.75pt}\put(-146,-155){
     \put(143,145)  {\scriptsize $ U^\vee $}
     \put(169,145)  {\scriptsize $ U $}
     }\setlength{\unitlength}{1pt}}
  \end{picture}}  
  \etb= d_U : U^\vee \oti U \rightarrow \one
  ~~,\qquad &
  \raisebox{-8pt}{
  \begin{picture}(26,22)
   \put(0,6){\scalebox{.75}{\includegraphics{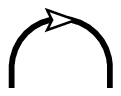}}}
   \put(0,6){
     \setlength{\unitlength}{.75pt}\put(-146,-155){
     \put(143,145)  {\scriptsize $ U $}
     \put(169,145)  {\scriptsize $ U^\vee $}
     }\setlength{\unitlength}{1pt}}
  \end{picture}}  
  \etb= \tilde d_U : U \oti U^\vee \rightarrow \one
  ~~,
\\[2em]
  \raisebox{-8pt}{
  \begin{picture}(26,22)
   \put(0,0){\scalebox{.75}{\includegraphics{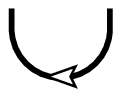}}}
   \put(0,0){
     \setlength{\unitlength}{.75pt}\put(-146,-155){
     \put(143,183)  {\scriptsize $ U $}
     \put(169,183)  {\scriptsize $ U^\vee $}
     }\setlength{\unitlength}{1pt}}
  \end{picture}}  
  \etb= b_U : \one \rightarrow U \oti U^\vee
  ~~,
  &
  \raisebox{-8pt}{
  \begin{picture}(26,22)
   \put(0,0){\scalebox{.75}{\includegraphics{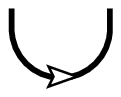}}}
   \put(0,0){
     \setlength{\unitlength}{.75pt}\put(-146,-155){
     \put(143,183)  {\scriptsize $ U^\vee $}
     \put(169,183)  {\scriptsize $ U $}
     }\setlength{\unitlength}{1pt}}
  \end{picture}}  
  \etb= \tilde b_U : \one \rightarrow U^\vee \oti U
  ~.
\eear
$$

A modular tensor category $\mathcal{C}$ is a semisimple abelian finite $\Cb$-linear ribbon category with a simple unit $\one$, $\text{End}(\one) =\Cb$ and satisfying an additional non-degeneration condition on braiding (given below). We denote the set of equivalence classes of simple objects in $\Cc$ by $I$, elements in $I$ by $i,j,k\in I$ and their representatives by $U_i, U_j, U_k$.  We define $\text{Dim}(\Cc):= \sum_i (\dim U_i)^2$. It is known that $\mathrm{Dim}\, \Cc \neq 0$ \cite{eno}. We also set $U_0=\one$. We define numbers $s_{i,j} \in \Cb$ 
by
$$  
 s_{i,j}   ~=~  \quad
\raisebox{-30pt}{
  \begin{picture}(110,65)
   \put(0,8){\scalebox{.75}{\includegraphics{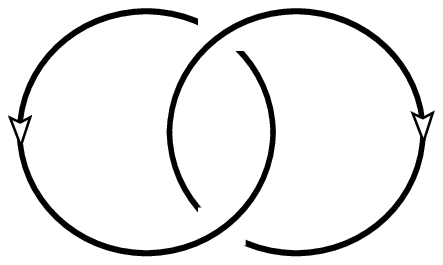}}}
   \put(0,8){
     \setlength{\unitlength}{.75pt}\put(-18,-19){
     \put( 98, 48)       {\scriptsize $ U_i $}
     \put( 50, 48)      {\scriptsize $ U_j $}
     }\setlength{\unitlength}{1pt}}
  \end{picture}}
~~.
$$
We also define $\dim(U) := \tilde{d}_U \circ b_U \in \Cb$. We have $s_{i,j}=s_{j,i}$ and $s_{0,i}=\dim U_i$. The non-degeneracy condition on the braiding of a modular tensor category is that the $|I|{\times}|I|$-matrix $s$ is invertible.

For $A \in \Cc$, we also choose a basis $\{ b_{A}^{(i;\alpha)} \}$ of $\Hom_{\Cc}(A, U_i)$
and the dual basis $\{ b_{(i;\beta)}^{A} \}$ of $\Hom_{\Cc}(U_i, A)$ for $i\in I$ such that
$b^{(i;\alpha)}_{A} \circ b_{(i;\beta)}^{A} = 
\delta_{\alpha\beta}\, \id_{U_i}$. We use the graphical notation
$$
b_A^{(i; \alpha)} = 
  \raisebox{-23pt}{
  \begin{picture}(30,52)
   \put(0,8){\scalebox{.75}{\includegraphics{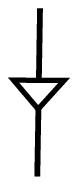}}}
   \put(0,8){
     \setlength{\unitlength}{.75pt}\put(-18,-11){
     \put(39,36)  {\scriptsize $ \alpha $}
     \put(23,65)  {\scriptsize $ U_i $}
     \put(23, 2)  {\scriptsize $ A $}
     }\setlength{\unitlength}{1pt}}
  \end{picture}}
\quad  , \qquad
b_{(i;\alpha)}^A =
  \raisebox{-23pt}{
  \begin{picture}(30,52)
   \put(0,8){\scalebox{.75}{\includegraphics{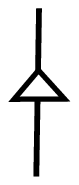}}}
   \put(0,8){
     \setlength{\unitlength}{.75pt}\put(-18,-11){
     \put(37,34)  {\scriptsize $ \alpha $}
     \put(20, 2)  {\scriptsize $ U_i $}
     \put(20,65)  {\scriptsize $ A $}
     }\setlength{\unitlength}{1pt}}
  \end{picture}} 
  ~~.
$$

Based on Huang and Lepowsky's earlier works on the tensor product of modules over a VOA \cite{hl-tcat-1}\cite{hl-tcat-2} \cite{hl-tcat-3}\cite{h-tcat-4}\cite{hl-tcat-5}, Huang proved the following theorem which is very important for us. 
\begin{thm}[\cite{h-mtc}]  \label{thm:huang-mtc}
If $V$ is a simple VOA satisfying
\bnu
\item $V_{(n)} =0$ for $n<0$, $V_{(0)}=\Cb \one$ and $V' \cong V$ as $V$-modules, 
\item Every $\Nb$-gradable weak $V$-module is completely reducible, 
\item $V$ is $C_2$-cofinite. 
\enu
Then we say $V$ is rational. Moreover, the category $\Cc_V$ of $V$-modules is a modular tensor category. 
\end{thm}

\medskip
\subsection{Non-holomorphic CFTs over $V^L\otimes V^R$}  \label{sec:cat-ffa}

An algebra in $\Cc$ or a $\Cc$-algebra is a triple $A=(A,m,\eta)$ where $A$ is an object of $\Cc$, $m$ (the multiplication) is a morphism $A \oti A \rightarrow A$ such that $m \cir (m \oti \id_A) \cir \alpha_{A,A,A} = m \cir (\id_A \oti m)$, and $\eta$ (the unit) is a morphism $\one \rightarrow A$ such that $m \cir (\id_A \oti \eta) = \id_A  \cir r_A$ and $m \cir (\eta\oti \id_A) = \id_A \cir l_A$. An algebra $A$ is called commutative if $m_A \circ c_{A,A} = m_A$. Similarly, one can define a coalgebra $A = (A,\Delta,\eps)$ where $\Delta : A \rightarrow A \oti A$ and $\eps : A \rightarrow \one$ obey a coassociativity and a counit conditions.

\begin{defn} {\rm 
A Frobenius algebra $A = (A,m,\eta,\Delta,\eps)$ is an algebra and a coalgebra such that the coproduct is an intertwiner of $A$-bimodules, i.e.\ 
$$(\id_A \oti m) \cir (\Delta \oti \id_A) = \Delta \oti m = (m \oti \id_A) \cir (\id_A \oti \Delta).$$ 
}
\end{defn}

We will use the following graphical representation for the morphisms of a Frobenius algebra,
$$
  m = \raisebox{-20pt}{
  \begin{picture}(30,45)
   \put(0,6){\scalebox{.75}{\includegraphics{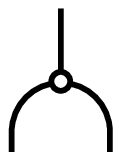}}}
   \put(0,6){
     \setlength{\unitlength}{.75pt}\put(-146,-155){
     \put(143,145)  {\scriptsize $ A $}
     \put(169,145)  {\scriptsize $ A $}
     \put(157,202)  {\scriptsize $ A $}
     }\setlength{\unitlength}{1pt}}
  \end{picture}}  
  ~~,\quad
  \eta = \raisebox{-15pt}{
  \begin{picture}(10,30)
   \put(0,6){\scalebox{.75}{\includegraphics{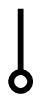}}}
   \put(0,6){
     \setlength{\unitlength}{.75pt}\put(-146,-155){
     \put(146,185)  {\scriptsize $ A $}
     }\setlength{\unitlength}{1pt}}
  \end{picture}}
  ~~,\quad
  \Delta = \raisebox{-20pt}{
  \begin{picture}(30,45)
   \put(0,6){\scalebox{.75}{\includegraphics{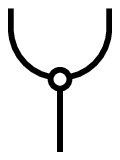}}}
   \put(0,6){
     \setlength{\unitlength}{.75pt}\put(-146,-155){
     \put(143,202)  {\scriptsize $ A $}
     \put(169,202)  {\scriptsize $ A $}
     \put(157,145)  {\scriptsize $ A $}
     }\setlength{\unitlength}{1pt}}
  \end{picture}}
  ~~,\quad
  \eps = \raisebox{-15pt}{
  \begin{picture}(10,30)
   \put(0,10){\scalebox{.75}{\includegraphics{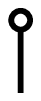}}}
   \put(0,10){
     \setlength{\unitlength}{.75pt}\put(-146,-155){
     \put(146,145)  {\scriptsize $ A $}
     }\setlength{\unitlength}{1pt}}
  \end{picture}}
  ~~.
$$
A Frobenius algebra $A$ in $\Cc$ is called symmetric if it satisfies the following identity: 
$$
\raisebox{-35pt}{
  \begin{picture}(50,75)
   \put(0,8){\scalebox{.75}{\includegraphics{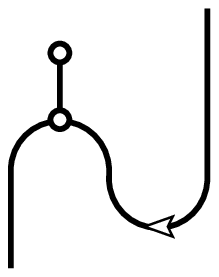}}}
   \put(0,8){
     \setlength{\unitlength}{.75pt}\put(-34,-37){
     \put(31, 28)  {\scriptsize $ A $}
     \put(87,117)  {\scriptsize $ A^\vee $}
     }\setlength{\unitlength}{1pt}}
  \end{picture}}  \quad
\,\, = \,\, \quad
  \raisebox{-35pt}{
  \begin{picture}(50,75)
   \put(0,8){\scalebox{.75}{\includegraphics{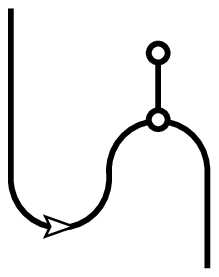}}}
   \put(0,8){
     \setlength{\unitlength}{.75pt}\put(-34,-37){
     \put(87, 28)  {\scriptsize $ A $}
     \put(31,117)  {\scriptsize $ A^\vee $}
     }\setlength{\unitlength}{1pt}}
  \end{picture}}
  ~~.
$$

Let $\Cc$ be a modular tensor category with simple objects $U_i, i \in I$ where $I$ is a finite set. If we replace the braiding and the twist in $\Cc$ by the antibraiding $c^{-1}$ and the antitwist $\theta^{-1}$ respectively, we obtain another ribbon category structure on $\Cc$. In order to distinguish these two distinct structures, we denote $(\Cc, c, \theta)$ and $(\Cc, c^{-1}, \theta^{-1})$ by $\Cc_+$ and $\Cc_-$ respectively. 

Let $\mathcal{D}$ be another modular tensor category. Let $W_j, j\in J$ be the representatives of simple objects in $\mathcal{D}$ where the index set $J$ is finite. By $\Cc \boxtimes \mathcal{D}$ we mean the tensor product of  $\Cb$-linear abelian categories \cite[def.\,1.1.15]{baki}, i.e.\ the category whose objects are direct sums of pairs $U \ti W$ of objects $U \in \Cc$ and $W \in \mathcal{D}$ and whose morphism spaces are
$$
 \Hom_{\Cc \boxtimes \mathcal{D}}(U \ti W,U' \ti W') = 
 \Hom_{\Cc}(U,U')\otimes_\Cb \Hom_{\mathcal{D}}(W,W')
$$
for pairs, and direct sums of these if the objects are direct sums of pairs. The category $\Cc \boxtimes \mathcal{D}_-$ is also a modular tensor category which has simple objects $U_i \times W_j$ for $i\in I$ and $j\in J$. In the case of $\Cc=\Cc_{V^L}$ and $\mathcal{D}=\Cc_{V^R}$, we have $U_i \times W_j = U_i\otimes_{\Cb} W_j$, where $\otimes_\Cb$ is the usual vector space tensor product.


\begin{thm}[\cite{ko-ffa}]  \label{thm:genus-0-cft}
The category of genus-0 CFTs (or self-dual full field algebras) over $V^L\otimes  V^R$ is isomorphic to the category of commutative symmetric Frobenius algebras in $(\Cc_{V^L})_+\boxtimes (\Cc_{V^R})_-$. 
\end{thm}

\begin{rema} \label{rema:cat-ffa} {\rm
The commutativity and the associativity of the Frobenius algebra in Theorem \ref{thm:genus-0-cft} is precisely the categorical formulation of the commutativity and the associativity of full field algebra discussed in Section \ref{sec:ffa}. 
}
\end{rema}

On a physical level of rigor, Sonoda \cite{sonoda} showed that higher genus theories only provide one additional condition: the modular invariance condition of 1-point genus-1 correlation functions \cite[Thm 3.8]{hk-modinv}.  For CFTs over $V^L\otimes V^R$, the categorical formulation of this {\it modular invariance condition} \cite{cardy-cond} is given in equation (\ref{eq:mod-inv}). Therefore, it is natural to have the following conjecture. 

\begin{conj}[\cite{cardy-cond}]
The category of CFTs over $V^L\otimes V^R$ is isomorphic to the category of 
commutative symmetric Frobenius algebras $A_\cl$ in $(\Cc_{V^L})_+\boxtimes (\Cc_{V^R})_-$ satisfying the following modular invariance condition: 
\be    \label{eq:mod-inv}
  \frac{\dim U_i \dim W_j}{\sqrt{\mathrm{Dim} \, \Cc_{V^L} \mathrm{Dim} \, \Cc_{V^R}}} ~~
\raisebox{-55pt}{
  \begin{picture}(90, 118)
   \put(0,8){\scalebox{.6}{\includegraphics{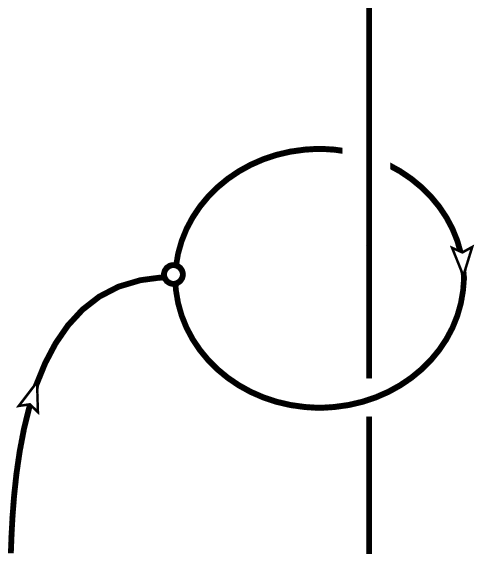}}}
   \put(0,8){
     \setlength{\unitlength}{.75pt}\put(-18,-19){
     \put(14, 8)     {\scriptsize $ A_\cl $}
     \put(126,80)  {\scriptsize $ A_\cl $}
     \put(95,8)    {\scriptsize $U_i \times W_j$ }
     \put(95,155){\scriptsize $U_i \times W_j $ }
     \put(62, 83)  {\scriptsize $m$} 
     }\setlength{\unitlength}{1pt}}
  \end{picture}}
\quad = \quad  \, \sum_{\alpha}
 \raisebox{-55pt}{
  \begin{picture}(90,118)
   \put(8,8){\scalebox{.6}{\includegraphics{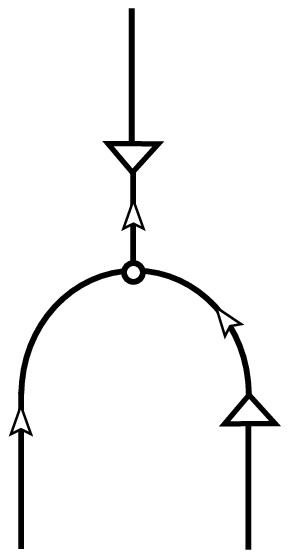}}}
   \put(0,8){
     \setlength{\unitlength}{.75pt}\put(-18,-19){
     \put(28, 8)  {\scriptsize $ A_\cl $}
     \put(50, 155)  {\scriptsize $ U_i \times W_j$}
     \put(80, 8)  {\scriptsize $ U_i \times W_j $}
     \put(36,95)  {\scriptsize $ A_\cl $}
     \put(85,75)    {\scriptsize $ A_\cl $}
     \put(67, 110)   {\scriptsize $\alpha$}
     \put(95, 50)     {\scriptsize $\alpha$}
     \put(53, 75)  {\scriptsize $m$}
}\setlength{\unitlength}{1pt}}
  \end{picture}}
\ee  
holds for all $i\in I$.
\end{conj}
 
\begin{rema} {\rm
Above categorical formulation of modular invariance condition is based on many earlier influential works on this subject, in particular Zhu's proof of the modular  transformation properties of $q$-traces on a rational VOA and its modules \cite{zhu}, Huang's proof of modular transformation properties of general genus-one correlation functions involving intertwining operators \cite{h-ode-torus}, Huang's proof of Verlinde conjecture \cite{h-ode-torus}\cite{h-verlinde} and Huang's construction of modular tensor categories \cite{h-mtc}. 
}
\end{rema}

The modular invariance condition is hard to check. In a special case it can be replaced by an easy-to-check condition. A $\Cc$-algebra $A$ is called haploid if $\dim \Hom(\one, A) =1$. 
\begin{thm}[\cite{cardy-alg}]
A haploid commutative symmetric Frobenius algebra $A_\cl$ in $\Cc_+\boxtimes \Cc_-$ is modular invariant if and only if 
$\dim A_\cl =\mathrm{Dim}(\Cc)$. 
\end{thm}

When $V^L=V^R=V$, examples of modular invariant commutative symmetric Frobenius algebras in $(\Cc_V)_+\boxtimes (\Cc_V)_-$ can be constructed explicitly. For example, let $U_i, i\in I$ be inequivalent irreducible $V$-modules. Then the object 
$$
A_\cl = \oplus_{i\in I} U_i^{\vee} \times U_i  \,\, \, \in (\Cc_V)_+\boxtimes (\Cc_V)_-,
$$
where $U_i^{\vee} \times U_i \cong U_i^{\vee} \otimes_\Cb U_i$ as vector spaces, 
has a structure of modular invariant commutative symmetric Frobenius $\Cc_{V\otimes V}$-algebra \cite{cardy-cond}\footnote{This result (but in a different framework) first appeared in \cite{tftcardycase} and later proved by Fjelstad, Fuchs, Runkel and Schweigert \cite{tft1,tft5} (see Remark \ref{rema:FRS}). The equivalence of two frameworks is known \cite{cardy-alg-2}.}. 
This construction is called the ``charge conjugation construction'' or the ``Cardy case''.  More constructions will be given in Section \ref{sec:examples}. In general, $V^L$ and $V^R$ can be different. For example, $V^R=V^L\otimes_\Cb V^{\natural}$ where $V^\natural$ is the Monster-Moonshine VOA.

\bigskip
\section{Open-closed CFTs}  \label{sec:op-cl-cft}

Open-closed (or boundary) CFT was first developed in physics by Cardy \cite{cardy:84,cardy:86,cardy:89ir,cardy:91tv}. It has important applications in the study of certain critical phenomena on surfaces with boundaries in condensed matter physics. It is also a powerful tool in the study of D-branes in string theory. Mathematically, it is easy to extend Segal's definition of a closed CFT to that of an open-closed CFT \cite{h-bcft}\cite{hukriz}. As we discussed in Section \ref{sec:cl-cft}, such a definition is not convenient if we want to apply the theory of VOA to the study of CFT. We will again take Vafa-Huang's approach to open-closed CFTs.

\subsection{Basic definitions} \label{sec:def-oc-cft}

\noindent {\bf The category $\mathcal{RS}_{\cl-\op}^p$}:  First, consider Riemann surfaces with unparametrized boundaries. A puncture living in the interior of the surface is called an interior puncture; that living on a boundary component is called a boundary puncture. The parametrization of an interior puncture is defined the same as before. The parametrization of a boundary puncture $p$ is the germ of a conformal map from neighborhoods of $p$ to neighborhoods of zero in the upper half-plane such that it is real analytic on the boundary, together with an orientation $\epsilon_p = \pm$. We introduce a unified notation for both interior punctures and boundary punctures: $(p, \sigma_p, f_p, \epsilon_p)$ where $\sigma_p$ takes the value of ``interior" or ``boundary", or equivalently ``closed string" and ``open string". A Riemann surface (not necessarily connected) with $k$ ordered positively oriented punctures (with the splitting $k=k^\cl+k^\op$ according to their $\sigma$-values) and $l$ ordered negatively oriented punctures (with $l=l^\cl + l^\op$) will be called a $(k^\cl, k^\op; l^\cl, l^\op)$-surface and will be denoted by
$$
\{ \Sigma | (p_1, \sigma_{p_1}, f_{p_1}, +), \dots, (p_k, \sigma_{p_k}, f_{p_k}, +) ;   (q_1, \sigma_{q_1}, f_{q_1}, -), \dots, (q_l, \sigma_{q_l}, f_{q_l}, -)\}
$$
or $\Sigma$ for simplicity. Two $(k^\cl, k^\op; l^\cl, l^\op)$-surfaces are conformally equivalent if there is a biholomorphic map between them that preserves the orders, the orientations and the local coordinates of the punctures. We denote the conformal equivalence classes of $\Sigma$ by $[\Sigma]$ and the moduli space of all $(k^\cl, k^\op; l^\cl, l^\op)$-surfaces by $\mathbb{E}_{(k^\cl, k^\op; l^\cl, l^\op)}$. We set $\mathbb{E}:=\cup_{k^\cl, k^\op; l^\cl, l^\op} \mathbb{E}_{(k^\cl, k^\op; l^\cl, l^\op)}$. On the one hand, the subspace of $\mathbb{E}$ containing only surfaces without boundary is nothing but $\mathcal{E}$. On the other hand, one can double a surface in $\mathbb{E}$ along its boundary components to obtain a closed surface with only interior punctures. This doubling defines a natural embedding of the moduli spaces: 
$$
\delta: \mathbb{E}_{(k^\cl, k^\op; l^\cl, l^\op)} \hookrightarrow \mathcal{E}_{(2k^\cl+k^\op, 2l^\cl+l^\op)}. 
$$
From now on, we identity $\mathbb{E}_{(k^\cl, k^\op; l^\cl, l^\op)}$ with its image in $\mathcal{E}_{(2k^\cl+k^\op, 2l^\cl+l^\op)}$. Notice that $\mathbb{E}$ inherits a smooth structure from that of $\mathcal{E}$. We also define 
$$
\mathbb{E}_{(k^\cl, k^\op; l^\cl, l^\op)}^g:= \mathbb{E}_{(k^\cl, k^\op; l^\cl, l^\op)} \cap \mathcal{E}_{(2k^\cl+k^\op, 2l^\cl+l^\op)}^g.
$$ 
We define $\mathbb{E}_{(k^\cl, k^\op; l^\cl, l^\op)}: = \cup_g \mathbb{E}_{(k^\cl, k^\op; l^\cl, l^\op)}^g$, $\mathbb{E}^g:=\cup_{k^\cl, k^\op; l^\cl, l^\op} \mathbb{E}_{(k^\cl, k^\op; l^\cl, l^\op)}^g$ and $\mathbb{E}:=\cup_g \mathbb{E}^g$.

The sewing operation defined on $\mathcal{E}$ is closed on the subspace $\mathbb{E}$ as long as we only allow sewing between oppositely oriented punctures with the same $\sigma$-values. Therefore, sewing operations are well-defined on $\mathbb{E}$. Now we are ready to give the following definition. 
\begin{defn}  {\rm 
The partial category $\mathcal{RS}_{\cl-\op}^p$ consists of
\bnu
\item objects: finite ordered sets $S$ with elements decorated by two colors: ``closed string" and ``open string", i.e. there is a splitting $S = S^\cl \cup S^\op$ such that $S^\cl \cap S^\op = \emptyset$. 

\item morphisms: $\Hom_{\mathcal{RS}_{\cl-\op}^p}(S_1, S_2) := \mathbb{E}_{(|S_1^\cl|, |S_1^\op|; |S_2^\cl|, |S_2^\op|)}$. 
\enu
The composition maps are defined by the sewing operations of Riemann surfaces.
}
\end{defn}

\begin{defn} {\rm 
An open-closed CFT is a real-analytic projective symmetric monoidal functor $\mathcal{F}: \mathcal{RS}_{\cl-\op}^p \to \mathcal{GV}$. 
}
\end{defn}

\subsection{Operads}

The category $\mathcal{RS}_{\cl-\op}^p$ is generated by two one-point sets: $\{ \mbox{a closed string} \}$, $\{ \mbox{an open string} \}$. We denote their $\Fc$-images by $V_\cl$, $V_\op$ respectively. 
We will be interested in a few substructures of $\mathcal{RS}_{\cl-\op}^p$:
\bnu

\item the set $D:= \cup_{n=0}^\infty\mathbb{E}_{0,n;0,1}^0$ gives a partial operad, 

\item the set $\mathbb{D}: = \cup_{m=0,n=0}^\infty \mathbb{E}_{0,n;0,m}^0$ gives a partial dioperad, 

\item the set $S:=\cup_{m=0,n=0}^\infty \mathbb{E}_{m,n;0,1}^0$ gives the so-called Swiss-cheese partial operad, 

\item the set $\mathbb{S}:=\cup_{k=0,l=0,m=0,n=0} \mathbb{E}_{k,l;m,n}^0$ gives the so-called Swiss-cheese partial dioperad. 

\enu

The structure $(V_\op, \{ \Fc([\Sigma]) \}_{[\Sigma]\in D})$\footnote{Only compositions associated to sewing of surfaces resulting in surfaces in $D$ are allowed.} is equivalent to a projective $D$-algebra structure on $V_\op$, and that of $(V_\op, \{ \Fc([\Sigma]) \}_{[\Sigma]\in \mathbb{D}})$ is equivalent to a projective $\mathbb{D}$-algebra. We will be interested in how to construct such algebras. The next two theorems tell us how to construct these algebras from two yet-to-be-introduced notions which will be discussed in the later subsections.

\begin{thm}[\cite{osvoa}]  \label{thm:osvoa}
An open-string VOA $A$ canonically gives a projective $D$-algebra. If $A$ is further equipped with a non-degenerate invariant bilinear form, then $A$ gives a projective $\mathbb{D}$-algebra. 
\end{thm}

\begin{rema}  {\rm 
The original theorem given in \cite{osvoa} says that the category of open-string VOAs of central charge $c$ is 
isomorphic to the category of holomorphic algebras (not projective!) over a certain partial-operad extension 
of $D$ by the $\frac{c}{2}$-th power of the determinant line bundle over $D$ and satisfying additional natural properties. We only give the result in above form to avoid technicalities.}
\end{rema}

\begin{thm}[\cite{ocfa}]   \label{thm:ocfa}
An open-closed field algebra $A$ canonically gives a projective $S$-algebra. If it is self-dual, then $A$ gives a projective $\mathbb{S}$-algebra.  
\end{thm}

\subsection{Open-string vertex operator algebras}  \label{sec:osvoa}

The notion of open-string VOA was introduced in \cite{osvoa}. 
We will not give a precise definition of it. Instead, we will use the structure of a projective $D$-algebra $(V_\op, \{ \Fc([\Sigma]) \}_{[\Sigma]\in D})$ to illustrate the basic ingredients and the properties of an open-string VOA. Note that there is a natural embedding $D \hookrightarrow K$ by doubling the disks in $D$ so that all the punctures are located on the equator of the resulting sphere. Therefore, we can adapt the same notation as (\ref{eq:element}) for elements in $D$. 

\bnu

\item An open-string VOA $V_\op$ is an $\mathbb{R}$-graded vector space such that $\dim (V_\op)_{(n)}< \infty$ and the grading is truncated from below, i.e. $(V_\op)_{(n)}=0$ for $n<<0$. 

\item Similar to VOA, we have $[\hat{\mathbb{C}}, (\infty, f_\infty, -))] \in D\subset K$ where $f_\infty: w\mapsto \frac{-1}{w}$. This gives rise to a distinguished element $\Fc([\hat{\mathbb{C}}, (\infty, f_\infty, -))])$ in $\overline{V}_\op$, denoted by $\one_\op$. In an open-string VOA, $\one_\op \in (V_\op)_{(0)}$. 

\item The $\Fc$ image of $[\hat{\mathbb{C}}, (0,f_0,+), (\infty, f_\infty, -))]$, where $f_0: w\mapsto w$ and $f_\infty: w\to \frac{-1}{w}$,  is $\id_{V_\op}$.

\item Consider the following element, 
$$
P(r):= [(\hat{\Cb}, (r, f_r, +), (0, f_0, +), (\infty, f_{\infty},-))] \in D \subset K
$$
where $r>0, f_r: w\rightarrow w-r, f_0: w\rightarrow w, f_{\infty}: w\rightarrow \frac{-1}{w}$.  The linear map $\mathcal{F}(P(z))$ gives a vertex operator 
$Y_\op(\cdot, r)\cdot := \mathcal{F}(P(r)): V_\op \otimes V_\op \rightarrow \overline{V}_\op$ 
for $r >0$. For $u,v\in V_\op$, we have (recall equation (\ref{eq:real-analytic}))
$$
u\otimes v \mapsto Y_\op(u,r)v = \sum_{n\in \Rb} u_n v \,\, r^{-n-1},
$$
where $u_n: V_\op \rightarrow \overline{V}$ and the nontrivial summands are those $n$ lying in $Q+\Zb$ where $Q$ is a finite set in $\Rb$. In an open-string VOA, $u_n\in \text{End}(V)$. Moreover, if $u\in (V_\op)_{(m)}$, then $u_n: (V_\op)_{(k)}\rightarrow V_{(k+m-n-1)}$. Therefore, $u_n v =0$ for $n>>0$. 

\item By the axioms of CFT, the following sum:
$$
\langle u', Y_\op(u_1, r_1)Y_\op(u_2, r_2)u_3\rangle := 
\sum_{n\in \Zb} \langle u', Y_\op(u_1, r_1)P_nY_\op(u_2, r_2)u_3\rangle
$$ 
is absolutely convergent in $\overline{V}_\op$ when $r_1>r_2>0$ for all $u_1,u_2,u_3\in V_\op$ and 
$u'\in V_\op'$. The following sum: 
$$
\hspace{1.5cm}\langle u', Y_\op(Y_\op(u_1, r_1-r_2)u_2, r_2) u_3\rangle := 
\sum_{n\in \Zb} \langle u', Y_\op(P_nY_\op(u_1, r_1-r_2)u_2, r_2) u_3\rangle
$$
is absolutely convergent in $\overline{V}_\op$ when $r_2>r_1-r_2>0$ all $u_1,u_2,u_3\in V_\op$ and 
$u'\in V_\op'$. 

\item We have the associativity of open-string VOA: 
$$
\langle u', Y_\op(u_1, r_1)Y_\op(u_2, r_2)u_3\rangle =
\langle u', Y_\op(Y_\op(u_1, r_1-r_2)u_2, r_2) u_3\rangle 
$$
when $r_1>r_2>r_1-r_2>0$ for all $u_1,u_2,u_3\in V$ and $u'\in V'$. This associativity is also called the OPE of boundary fields in physics. 

\item Consider the element $\Sigma_{\omega}:=[(\hat{\Cb}, (\infty, f_{\infty}, -))]\in D$ with $f_{\infty}: w\rightarrow \frac{-1}{w-\eps}$. It is clear that $\mathcal{F}(\Sigma_{\omega})$ is an element in $\overline{V}_\op$. Then 
$\omega_\op:=\frac{d}{d\eps}|_{\eps=0} \mathcal{F}(\Sigma_{\omega})$ gives a distinguished element in $\overline{V}_\op$. Actually, $\omega_\op\in (V_\op)_{(2)}$. One has
$$
Y_\op(\omega_\op, r) = \sum_{n\in \Zb} L(n) r^{-n-2}~,
$$
where $L(n), n\in \Zb$, can be shown to generate a Virasoro Lie algebra, i.e. 
$$
\hspace{1cm} [L(m), L(n)] = (m-n)L(m+n) + \frac{c}{12}(m^3-m) \delta_{m+n, 0}~, \quad\quad \forall m,n\in \Zb,
$$
where $c\in \Cb$ is called the central charge. Moreover, $\omega=L(-2)\one$ and $L(0)$ gives the grading operator. 

\item We thus denote an open-string VOA by its basic ingredients: $(V_\op, Y_\op, \one_\op, \omega_\op)$. 

\enu

It is also straight forward to define the notion of a module over an open-string VOA. We will not do it here. Instead, we will explore a relation between the representation theory of VOAs and open-string VOAs. As we will show later that an open-string VOA is always a module over a VOA. This fact is extremely useful when the underlying VOA is rational. 

We introduce a formal vertex operator $\Y_\op$ such that 
$$
\Y_\op(u, x)v = \sum_{n\in \Rb} u_n v x^{-n-1}.
$$
for $u, v\in V_\op$. Namely, $\Y_\op(u, x)v|_{x=r} = Y_\op(u, r)v$. We define the so-called {\it meromorphic center} $C_0(V_\op)$ of the open-string VOA $V_\op$ as follows: 
\bea
C_0(V_\op) &:=& \Big\{ u \in \coprod_{n\in \Zb} V_{(n)} | \Y_\op(u, x) \in (\mathrm{End} V) [[x, x^-1]],    \nn
&& \hspace{1.5cm}  \Y_\op(v, x) u = e^{xL(-1)} \Y_\op(u, -x) v, \forall v\in V_\op \Big\}   \nonumber
\eea
It is easy to show that $\one_\op, \omega_\op \in C_0(V_\op)$ and the image of $C_0(V_\op)\otimes C_0(V_\op)$ under $\Y_\op|_{C_0(V_\op)}$ lies in $C_0(V_\op)[[x, x^{-1}]]$. Moreover, we have the following result: 
\begin{prop}
The quadruple $(C_0(V_\op), \Y_\op|_{C_0(V_\op)}, \one_\op, \omega_\op)$ is a VOA. $V_\op$ is a module over $C_0(V_\op)$ and 
$\Y_\op$ is an intertwining operator of type $\binom{V_\op}{V_\op V_\op}$ \cite{FHL} with respect to $C_0(V_\op)$. 
\end{prop}

By the above results, an open-string VOA can be constructed and studied by the representation theory of VOAs. An open-string VOA containing a VOA $V$ as a sub-VOA of its meromorphic center is called an open-string VOA over $V$. We denote the canonical embedding $V \hookrightarrow V_\op$ by $\iota_\op$. If $V$ is a rational VOA, one can apply the tensor category theory to construct open-string VOAs. 

\medskip
An invariant bilinear form $(\cdot, \cdot)_\op: V_\op \otimes V_\op \to \Cb$ is called invariant if the following identity: 
$$
(w, Y_\op(u, z)v)_\op = ( e^{-r^{-1}L(-1)} Y_\op(w, r^{-1})e^{-rL(1)} r^{-2L(0)}u, v)_\op
$$
holds for $u, v, w\in V_\op$ and $r>0$. An open-string VOA equipped with a non-degenerate symmetric invariant bilinear form is called self-dual. It was shown in \cite{cardy-cond} that a self-dual open-string VOA canonically produces a projective algebra over the partial dioperad $\mathbb{D}$. 

\subsection{Open-closed field algebras}  \label{sec:ocfa}

The notion of open-closed field algebra was introduced in \cite{ocfa}. We will not give an explicit definition of it. Instead, by Theorem \ref{thm:ocfa}, we will illustrate the structures and the properties of an open-closed field algebra via the structures on $S$. 
\bnu

\item The Swiss-cheese partial operad $S$ contains the spherical partial operad $K$ as a subset. Therefore, an open-closed field algebra contains a full field algebra $V_\cl$ as a substructure. 

\item The Swiss-cheese partial operad $S$ also contains the partial operad $D$ as a subset. Therefore, an open-closed field algebra contains an open-string VOA $V_\op$ as a substructure. 

\item Consider the following element, for $z\in \mathbb{H}$, 
$$
D(z):=[(\Sigma|(z, f_z, +), (0,f_0,+), (\infty, f_\infty, -))],
$$
where $f_z: w \mapsto w-z$, $f_0: w\mapsto w$, $f_\infty: w\mapsto \frac{-1}{w}$. It is clear that
$\Fc(\Sigma): V_\cl \otimes V_\op \to \overline{V}_\op$. We set $Y_\co(\cdot; z, \bar{z}) \cdot :=\Fc(\Sigma)$. 

\item Since the entire Swiss-cheese partial operad $S$ is generated by the elements in $K$, $D$ and $D(z)$, an open-closed field algebra can be described by the data $V_\op, V_\cl, Y_\co$. Therefore, we denote an open-closed field algebra by a triple $(V_\op|V_\cl, Y_\co)$. 

\enu

The notion of open-closed field algebra introduced in \cite{ocfa} is too general for the purpose of this paper. We will only consider the so-called {\it analytic open-closed field algebras} which satisfy the following analytic properties:  
\bnu

\item 
\bnu

\item $Y_\op$ can be extended to a map $V_\op \otimes V_\op \times (\Cb^\times/-\Rb_+) \to \overline{V}_\op$:
$$
Y_\op(v, r) = Y_\op(u, z) |_{z=r}; 
$$

\item $Y_\co$ can be extended to a map
$V_\cl\otimes V_\op \times \Hb \times \overline{\Hb} \rightarrow
\overline{V}_\op$ such that for $z\in \Hb, \zeta\in \overline{\Hb}$
$$
Y_\co(u; z, \bar{z}) = Y_\co(u;z,\zeta) 
|_{\zeta=\bar{z}}. 
$$
\enu

\item For $n\in \Nb$, $v_1,\cdots, v_{n+1}\in V_\op, v'\in (V_\op)'$ 
and $u_1,\cdots, u_n\in V_\cl$ and $z_1, \cdots, z_n\in \Hb, \zeta_1, \cdots, \zeta_n\in \overline{\Hb}$, the series
$$
\langle v', Y_\co(u_1; z_1,\zeta_1)Y_\op(v_1, r_1) \cdots 
Y_\co(u_n; z_n,\zeta_n)Y_\op(v_n, r_n)v_{n+1}\rangle
$$
is absolutely convergent when 
$|z_1|, |\zeta_1| > r_1> \cdots >|z_n|, |\zeta_n|>r_n>0$ and 
can be extended to a (possibly multivalued) analytic function on:
$$
\{ (z_1, \zeta_1, r_1, \dots, z_n, \zeta_n, r_n) \in M_{\Cb}^{3n} \}
$$ 
where $M_{\Cb}^{n} := \{ (z_1, \dots, z_n) \in \Cb^n | z_i\neq z_j, \,  \,
\mbox{for $i,j=1, \dots, n$ and $i\neq j$} \}$.

\item For $n\in \Nb$, $v', u_1, \dots, u_{n+1}\in V_\cl$ and $z_1, \cdots, z_n, \zeta_1, \cdots, \zeta_n\in \Cb$
the series
$$
\langle v', Y_\cl(u_1; z_1, \zeta_1)\cdots 
Y_\cl(u_n; z_n, \zeta_n)u_{n+1}\rangle 
$$
is absolutely convergent when 
$|z_1|>\dots >|z_{n}|>0$ and $|\zeta_1|>\dots >|\zeta_n|>0$
and can be extended to an analytic function on $M_{\Cb}^{2n}$.

\item For $v'\in V_\op', v_1, v_2\in V_\op, u\in V_\cl$ and $z\in \Hb, \zeta\in \overline{\Hb}$, the series
$$
\langle v', Y_\op(Y_\co(u; z, \zeta)v_1, r)v_2\rangle
$$
is absolutely convergent when $r>|z|, |\zeta|>0$.

\item For $v'\in V_\op', v\in V_\op, u_1, u_2\in V_\cl$ and $z\in \Hb, \zeta\in \overline{\Hb}$, the series
$$
\langle v', Y_\co(Y_\cl(u_1; z_1, \zeta_1)u_1; 
z_2, \zeta_2)v\rangle,
$$
converges absolutely when $|z_2|>|z_1|>0$, $|\zeta_2|>|\zeta_1|>0$
and $|z_1|+|\zeta_1|<|z_2-\zeta_2|$. 

\enu

An analytic open-closed field algebra satisfies the following nice properties: 

\begin{enumerate}

\item Unit property: $Y_\co(\one_\cl; z, \zeta) = \id_{V_\op}$. 

\item Associativity I:  
For $u\in V_\cl, v_1,v_2\in V_\op, v'\in V_\op'$ and
$z\in \Hb,\zeta\in \overline{\Hb}$,  we have
$$
\langle v', Y_\co(u; z, \zeta)
Y_\op(v_1, r)v_2\rangle  
=\langle v', Y_\op
( Y_\co(u; z-r, \zeta-r)v_1, r)v_2\rangle 
$$
when $|z|,|\zeta|>r>0$ and $r>|r-z|,|r-\zeta|>0$. 

\item Associativity II: 
For $u_1, u_2 \in V_\cl, v_1,v_2\in V_\op, v'\in V_\op'$ and
$z_1,z_2\in \Hb,\zeta_1, \zeta_2\in \overline{\Hb}$, we have 
\bea  \label{asso-co-co}
&& \langle w', Y_\co (u_1; z_1, \zeta_1)
Y_\co(u_2; z_2, \zeta_2)v_2\rangle  \nn
&& \hspace{1cm} =\langle v', Y_\co
( Y_\cl(u_1; z_1-z_2, \zeta_1-\zeta_2)
u_2; z_2, \zeta_2)v_2\rangle    \nonumber
\eea
when $|z_1|, |\zeta_1|> |z_2|, |\zeta_2|$ and 
$|z_2| > |z_1-z_2|>0, |\zeta_2|>|\zeta_1-\zeta_2|>0$ and
$|z_2-\zeta_2|>|z_1-z_2|+|\zeta_1-\zeta_2|$.

\item Commutativity I: The map $Y_\co$ can be uniquely extended to 
$V_\cl\otimes V_\op\times R$, where 
$$
R : = \{ (z, \zeta) \in \Cb^2 | z\in \Hb \cup \Rb_+, \zeta\in \overline{\Hb}\cup \Rb_+, z\neq \zeta \}. 
$$
For $u\in V_\cl, v_1,v_2\in V_\op$ and $v'\in V_\op'$,  
$$ 
\langle v', Y_\co(u; z, \zeta)
Y_\op(v_1, r)v_2\rangle,
$$
which is absolutely convergent when $z>\zeta > r>0$, and 
$$
\langle v', Y_\op(v_1, r)Y_\co(u; z, \zeta)v_2\rangle,
$$
which is absolutely convergent when $r>z>\zeta>0$, are the
analytic continuations of each other along the following path 
\footnote{The extended domain $R\backslash \{z=r, \text{or}, \zeta =r\}$ 
is simply connected for fixed $r>0$, 
all possible paths of analytic continuation  
are homotopically equivalent.}. 
$$  
\epsfxsize 0.4\textwidth
\epsfysize 0.15\textwidth
\epsfbox{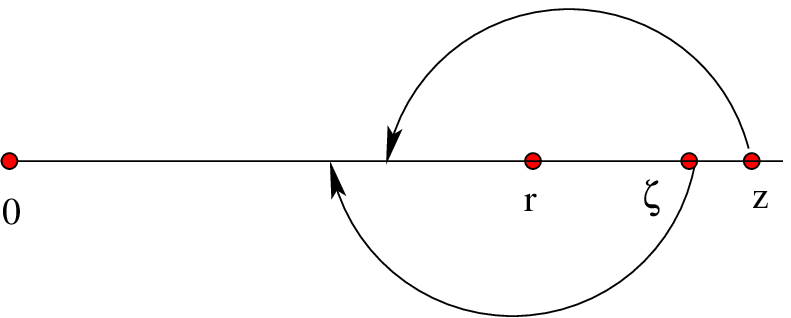}
$$

\item Commutativity II: For $u_1,u_2\in V_\cl, v\in V_\op$ and $v'\in V_\op'$,
$$  
\langle v', Y_\co(u_1; z_1, \zeta_1)
Y_\co(u_2; z_2, \zeta_2)v\rangle,
$$
which is absolutely convergent when $z_1>\zeta_1>z_2>\zeta_2>0$, and
$$
\langle v', Y_\co(u_2; z_2, \zeta_2)
Y_\co(u_1; z_1, \zeta_1)v\rangle, 
$$
which is absolutely convergent when $z_2>\zeta_2>z_1>\zeta_1>0$, 
are analytic continuation of each other along the following paths.
$$ 
\hspace{1cm}\epsfxsize 0.8\textwidth
\epsfysize 0.15\textwidth
\epsfbox{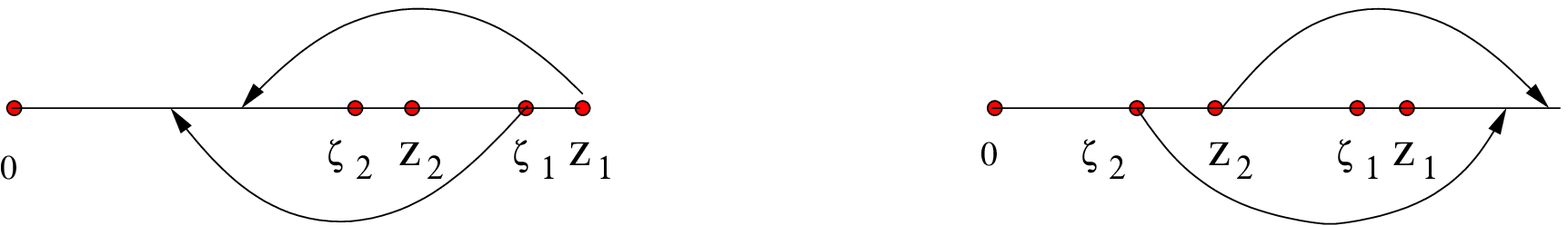}
$$

\end{enumerate}

\begin{rema}  {\rm
The commutativity II follows from the Associativity II and the commutativity of the full field algebra. This is an analogue of the statement that the action of a commutative algebra on its module is automatically commutative. 
}
\end{rema}

For a given open-closed field algebra $(V_\op|V_\cl, Y_\co)$, the data $Y_\co$ can also be described by slightly simpler data. We define $\iota_\co: V_\cl \times \Hb \times \overline{\Hb} \rightarrow \overline{V}_\op$ as follows: 
$$
\iota_\co^{(z,\zeta)}(u) := Y_\co(u; z, \zeta) \one_\op. 
$$
for $u\in V_\cl$ and $v\in V_\op$. Conversely, $Y_\co(u; z, \zeta)$ can be obtained from $\iota_\co$ by defining 
$$
Y_\co(u; z, \zeta):=Y_\op(\iota_\co^{(z-r, \zeta-r)}(u), r)v
$$ 
for $r>|z-r|,|\zeta-r| >0$ and by the analytic continuations. Therefore, we can also denote the open-closed field algebra by a triple $(V_\op|V_\cl, \iota_\co)$. Moreover, it is easy to show that this map $\iota_\co$ plays the role of an algebra homomorphism between $V_\cl$ and $V_\op$. It becomes evident in the classification of open-closed rational CFTs (see Theorem \ref{thm:cardy-alg}).

\begin{rema}  \label{def:chiral-module}  {\rm
In an open-closed field algebra $(V_\op|V_\cl, \iota_\co)$, if we ignore the algebra structure on $V_\op$, just viewing it as a module over one copy of the Virasoro algebra equipped with an action $Y_\co$ of $V_\cl$ satisfying the Associativity II and the unit property, we obtain some kind of a ``$V_\cl$-module" structure on $V_\op$ given by $Y_\co$. We will call it a chiral $V_\cl$-module. Due to the fact that the map $\iota_\co$ is some kind of an algebraic homomorphism, a $V_\op$-module is automatically a chiral $V_\cl$-module. 
}
\end{rema}

\medskip
From now on, we will be interested in open-closed field algebras over a VOA $V$, which consists of a full field algebra over $V\otimes V$ $(V_\cl, V\otimes V \xrightarrow{\iota_\cl} V_\cl)$ and an open-string VOA over $V$ $(V_\op, V\xrightarrow{\iota_\op} V_\op)$. In this case, one can define two maps $h^L : V \rightarrow V_\op$ and $h^R:  V \rightarrow V_\op$ as follows: for all $u, v \in V$, 
\bea  \label{h-L-R}
&&h^L: u \mapsto \lim_{z\rightarrow 0}
Y_\co(\iota_\cl(u\otimes \one), z)\one_\op, \nn
&&h^R: v \mapsto \lim_{\zeta\rightarrow 0}
Y_\co(\iota_\cl(\one \otimes v), \zeta)\one_\op.  \nonumber
\eea
Notice that $h^L, h^R$ preserve the conformal weights. Namely
\bea
\wt \, h^L(u) &=& \wt \, (u_{-1}\, \one_\op) = \lwt \, u,    \nn
\wt \, h^R(v) &=& \wt \, (v_{-1}\, \one_\op) = \rwt \, v.   \nonumber
\eea
Therefore, both $h^L$ and $h^R$ can be naturally extended
to maps $\overline{V} \rightarrow \overline{V_\op}$.
We still denote the extended maps by $h^L$ and $h^R$ respectively.

\begin{defn} \label{def-ocfa-V}
{\rm 
Let $(V, Y, \one, \omega)$ be a vertex operator algebra.
An {\it open-closed field algebra over $V$} is an 
analytic open-closed field algebra 
$$
( V_\op|V_\cl, Y_\co),
$$ 
where $V_\cl=(V_\cl, Y_\cl, \one_\cl, \omega_\cl)$ 
is a full field algebra over 
$V\otimes V$ and $(V_\op, Y_\op, \one_\op, \omega_\op)$ 
is an open-string vertex operator algebra
over $V$, satisfying the following conditions:
\bnu
\item  {\em $V$-invariant boundary condition}: 
$h^L=h^R=\iota_\op$.

\item {\em Chirality splitting property}: 
$\forall u\in V_\cl$,  
$u=u^L\otimes u^R\in W^L\otimes W^R \subset V_\cl$ 
for some $V$-modules $W^L, W^R$. 
There exist $V$-modules $W_1, W_2$ and intertwining operators 
$\Y^{(1)}, \Y^{(2)}, \Y^{(3)}, \Y^{(4)} $ of type 
$\binom{V_\op}{W^LW_1}$, 
$\binom{W_1}{W^RV_\op}$, $\binom{V_\op}{W^RW_2}$,
$\binom{W_2}{W^LV_\op}$ respectively, such that
$$
\langle w', Y_\co(u; z, \zeta)w\rangle
=  \langle w', \Y^{(1)}(u^L, z)\Y^{(2)}(u^R, \zeta) w\rangle
$$
when $|z|>|\zeta|>0$, 
and 
$$
\langle w', Y_\co(u; z, \zeta)w\rangle
= \langle w', \Y^{(3)}(u^R, \zeta)\Y^{(4)}(u^L, z) w\rangle
$$
when $|\zeta|>|z|>0$
for all $u\in V_\cl, w\in V_\op, w'\in V_\op$.  
\enu    

In the case that $V$ is generated by $\omega$, i.e. $V=\langle \omega \rangle$, the $\langle \omega\rangle$-invariant boundary condition is also called the {\it conformally invariant boundary condition}. We also call an open-closed field algebra over $\langle \omega\rangle$ an {\it open-closed conformal field algebra}. 
}
\end{defn}

\begin{rema} {\rm 
In order to obtain an open-closed CFT, the minimal requirement on the boundary is the conformally invariant boundary condition. By restricting ourselves to $V$-invariant boundary conditions and assuming $V$ is rational, we will be able to obtain a nice categorical formulation of open-closed field algebra over $V$. The price we pay is that there are too few $V$-invariant D-branes to reveal an interesting geometry. To classify all conformally invariant D-branes for a given closed rational CFT is a central task in the program of SAG. 
}
\end{rema}

\begin{defn}  {\rm
An open-closed field algebra $(V_\op|V_\cl, Y_\co)$ is called self-dual if there are non-degenerate symmetric invariant bilinear forms $(\cdot, \cdot)_\op$ and $(\cdot, \cdot)_\cl$ on $V_\op$ and $V_\cl$ respectively. 
}
\end{defn}

For a self-dual open-closed field algebra $(V_\op|V_\cl, \iota_\co)$, it is easy to obtain the Ishibashi states and the boundary states. Indeed, for $v\in V_\op$, we define the boundary state $\mathbf{B}(v) \in \overline{V_\cl}$ associated to $v$ by 
$$
\mathbf{B}(v) = e^{L(-1)} (\bar{z}_0-z_0)^{L(0)} \otimes
e^{L(-1)} \overline{\bar{z}_0-z_0}^{\, L(0)}
\iota_\co^*(z_0, \bar{z}_0) (v),
$$
where $z_0 \in \Hb$ and $\iota_\co^\ast$ is the adjoint of $\iota_\co$, i.e. for $v\in V_\op$ and $u\in V_\cl$, 
$$
(v, \iota_\co^{(z, \bar{z})} (u))_\op = ((\iota_\co^{(z, \bar{z})})^\ast(v), u)_\cl. 
$$
Actually, $\mathbf{B}(v)$ is independent of $z_0$. This explains the notation. 
\begin{prop}[Ishibashi States \cite{cardy-cond}]
If $v\in V_\op$ is such that $L(-1)v=0$, then, for $z_0\in \Hb$,  $\mathbf{B}(v)$ is an Ishibashi state, i.e.
$$
(L^L(n) - L^R(-n)) \mathbf{B}(v) = 0, \quad\quad \forall n\in \Zb. 
$$
\end{prop}

An open-closed CFT containing only spheres and disks with arbitrary number of incoming and outgoing interior and boundary punctures is equivalent to a projective algebra over the Swiss-cheese partial dioperad \cite{ocfa}. Since spheres and disks (under the doubling map) cover all genus-zero surfaces, a projective algebra over the Swiss-cheese partial dioperad is also called a genus-zero open-closed CFT. 
\begin{thm}[\cite{ocfa}]
A self-dual open-closed field algebra over a rational VOA $V$ canonically gives a genus-zero open-closed CFT.  
\end{thm}

Notice that all the surfaces in open-closed CFT can be obtained by gluing spheres and disks. Therefore, an open-closed field algebra over a rational $V$ contains all the building blocks of an open-closed CFT. Higher genus surfaces (in the sense of the doubling map) only provide certain compatibility conditions. The modular invariance condition discussed in Section \ref{sec:cl-cft} is one of such compatibility conditions. It has long been conjectured that the only remaining compatibility condition is the so-called Cardy condition \cite{Lew}. The Cardy condition comes from two ways of realizing the same world-sheet: 
\vspace{2mm}
\begin{center}
\includegraphics[width=0.85\textwidth]{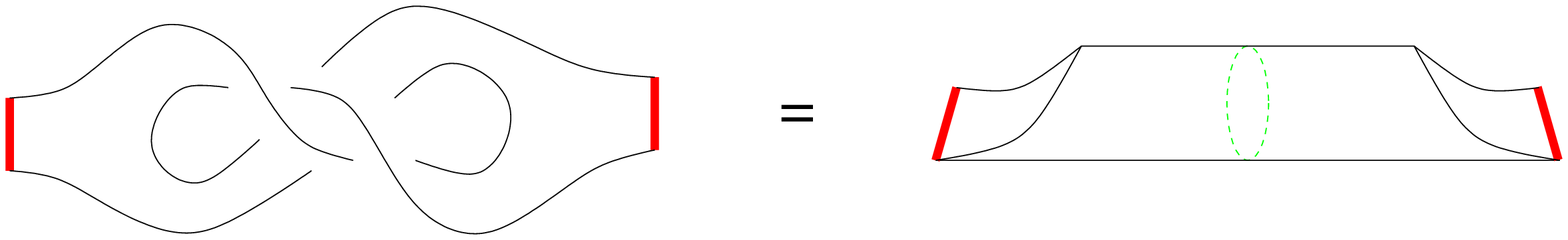}
\end{center}
\vspace{2mm}
in which the two (red) line segments at two sides of the above surfaces represent incoming or outgoing open strings. A precise definition of the Cardy condition in terms of correlation functions has been worked out in \cite[Def 3.4, Thm 3.10]{cardy-cond}. We will not recall it here. We only need its categorical formulation \cite[eq. (5.31)]{cardy-cond}\cite[eq. (3.14)]{cardy-alg}, which will be given in the next section.  

\begin{conj}  \label{conj:oc-alg-bcft}
An open-closed field algebra over a rational VOA $V$ satisfying the modular invariance condition and the Cardy condition gives a consistent open-closed CFT. 
\end{conj}

\subsection{Classification of open-closed CFTs over $V$}

Now we assume that $V$ is rational. Namely, the category $\Cc_V$ of $V$-modules is a modular tensor category. In this case, $\Cc_V \boxtimes (\Cc_V)_-$ is canonically equivalent to the monoidal center $Z(\Cc_V)$ of $\Cc_V$ as modular tensor categories. For simplicity, we set $\CcC:=\Cc_V \boxtimes (\Cc_V)_-$. We denote the tensor product functor $A \times B \to A\otimes B$ from $\CcC$ to $\Cc_V$ by $T$. We denote the right adjoint of $T$ by $R$. Actually, $R$ is also a left adjoint of $T$. $T$ is a monoidal functor. It was shown in \cite{cardy-alg} that $R$ is both a lax and a colax functor. Moreover, it is a Frobenius functor such that it maps a Frobenius algebra in $\Cc_V$ to a Frobenius algebra in $\Cc_V \boxtimes (\Cc_V)_-$. 

We need the notion of center of an algebra $A$ in $\Cc_V$ in order to state our classification result. We first define what a left center of an algebra in $\CcC$. We denote the braiding of $\CcC$ by $c_{X,Y}: X\otimes Y \to Y\otimes X$. Then the left center of an algebra $(B, m_B, \iota_B)$ in $\CcC$ is the maximal subobject 
$e: Z \hookrightarrow B$ such that $m_B \circ (e \otimes \id_B) = m_B \circ c_{B, B} \circ (e \otimes \id_B)$. 
\begin{defn}[\cite{unique}] {\rm
The center of $A$ in $\Cc_V$ is the left center $C_l(R(A))$ of $R(A)$ in $\CcC$}.
\end{defn}

\begin{rema} {\rm
The above definition of center of an algebra is not obviously natural. We will briefly describe its naturalness in this remark via its universal property. The monoidal center $Z(\Cc)$ of a tensor category $\Cc$ is the category of pairs $(Z, z)$ where $Z\in \Cc$ and $z: Z\otimes - \to -\otimes Z$ is the half braiding (satisfying certain properties). We denote the forgetful functor $Z(\Cc) \to \Cc$ by $F$. When $\Cc$ is a modular tensor category, there is a canonical equivalence of ribbon categories $\phi: \CcC \xrightarrow{\cong} Z(\Cc)$. Moreover, we have $T \cong F \circ \phi$. Therefore, one can equivalently define the center of $A$ in $\Cc$ as an object in $Z(\Cc)$. More precisely, Davydov \cite{davydov} showed that the center of $A$ can be equivalently defined by a pair $((Z, z), Z \xrightarrow{e_Z} A)$, where $(Z, z)$ is an object in $Z(\Cc)$ and $e_Z$ is a morphism in $\Cc$, such that it is terminal among all such pairs $( (X,x), e_X)$ satisfying the following commutative diagram: 
$$
\xymatrix{
X \otimes A \ar[r]^{e_X\otimes \id_A}  \ar[d]_{x} & A \otimes A \ar[dr]^{m_A}  &  \\
A \otimes X \ar[r]^{\id_A \otimes e_X} & A\otimes A \ar[r]^{m_A} & A \, .
}
$$
Similar to the case of algebras over vector spaces, the above definition of center is equivalent to a certain internal endomorphism of the identity functor on the category $A$-modules \cite{davydov}. 
}
\end{rema}

\medskip
Given a morphism $f: A\to B$ between two Frobenius algebras $(A, m_A, \eta_A, \Delta_A, \epsilon_A)$ and $(B, m_B, \eta_B, \Delta_B, \epsilon_B)$, we define the right adjoint of $f$ to be
$$
f^\ast :=  ( (\epsilon_B \circ m_B) \otimes \id_A) \circ (\id_B \otimes f \otimes \id_A) \circ (\id_B \otimes (\Delta_A \circ \eta_A)). 
$$

Now we are ready to give a categorical formulation of open-closed field algebra over a rational $V$ satisfying additional properties. This result was first obtained in \cite[Thm. 5.15]{cardy-cond}. The following version is taken from \cite{cardy-alg}. 
\begin{thm}[\cite{cardy-cond}\cite{cardy-alg}] \label{thm:cardy-alg}
An open-closed field algebra over $V$ satisfying the modular invariance condition and the Cardy condition  is equivalent to a triple
$$
(A_\op|A_\cl, \iota_\co)
$$
where
\bnu
\item $A_\op$ is a symmetric Frobenius algebra in $\Cc_V$; 
\item $A_\cl$ is a modular invariant commutative symmetric Frobenius algebra in $\CcC$; 
\item $\iota_\co: A_\cl \rightarrow R(A_\op)$ is an algebra homomorphism factoring through $e: Z(A_\op)\hookrightarrow R(A_\op)$
\enu
satisfying the Cardy condition: 
\be   \label{eq:cardy-CC}
  \iota_\co \circ \iota_\co^\ast ~=~
  \raisebox{-40pt}{
  \begin{picture}(54,80)
   \put(0,8){\scalebox{.75}{\includegraphics{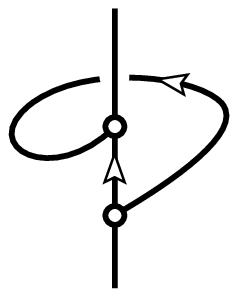}}}
   \put(0,8){
     \setlength{\unitlength}{.75pt}\put(-18,-19){
     \put(45, 10)  {\scriptsize $ R(A_\op) $}
     \put(45,105)  {\scriptsize $ R(A_\op) $}
     \put(80, 80)  {\scriptsize $ R(A_\op) $}
     }\setlength{\unitlength}{1pt}}
  \end{picture}}
~~.
\ee
where $\iota^*$ is the right adjoint of $\iota_\co$. 
\end{thm}

Following \cite{cardy-cond}\cite{cardy-alg}, we will call the triple defined in Theorem \ref{thm:cardy-alg} a Cardy $\Cc_V|\CcC$-algebra. 

\begin{rema} {\rm 
If Conjecture \ref{conj:oc-alg-bcft} is correct, Theorem \ref{thm:cardy-alg} gives a classification of open-closed rational CFTs satisfying the $V$-invariant boundary condition. A similar classification result has been obtained in another approach towards rational CFTs via the theory of conformal nets \cite{longo-rehren-1,longo-rehren-2, rehren}. In this approach, no higher genus surface is involved. The modular invariance condition and the Cardy condition is somehow captured by a so-called Haag dual condition (see \cite{kr-cft} for a review).  
}
\end{rema}

\begin{defn} {\rm
For a given full field algebra $A_\cl$ over $V\otimes V$, a {\it $V$-invariant D-brane} is a symmetric Frobenius algebra $A$ in $\Cc_V$ together with an algebra homomorphism $\iota_\co: A_\cl \to R(A)$ such that the triple $(A|A_\cl, \iota_\co)$ gives a Cardy $\Cc_V|\CcC$-algebra. 
}
\end{defn}

It is easy to see that all $V$-invariant D-branes naturally form a category instead of a set. In general, such $V$-invariant D-branes are too few to give an interesting geometry. For later purposes, one should consider symmetry-broken D-branes, i.e. $U$-invariant D-branes where $U$ is a sub-VOA of $V$. In particular, when $U=\langle \omega \rangle$, they are also called {\it conformally invariant D-branes}. Therefore, to have a good picture of SAG, we need to classify all the conformally invariant D-branes associated to a closed CFT. In particular, we need to understand how $V$ and $V$-modules decompose into modules over $\langle \omega \rangle$ the Virasoro sub-VOA of $V$. But in general, $\langle \omega \rangle$ is very small and irrational. Therefore, in order to understand SAG associated to a closed rational CFT, it is still unavoidable to study irrational theory.

\subsection{Constructions}  \label{sec:examples}

An algebra $A$ in a modular tensor category $(\Cc, \otimes, \one)$ is called haploid if 
$\dim \Hom_\Cc(\one, A) =1$. 
\begin{defn}   \label{def:separable} {\rm
A Frobenius algebra $(A, m_A, \eta_A, \Delta_A, \epsilon_A)$ is called separable if 
$$
m_A \circ \Delta_A = \beta_A \id_A \quad \mbox{and} \quad
\epsilon_A \circ \eta_A = \beta_A' \id_\one \quad
\mbox{and}\quad \beta_A, \beta_A' \in \Cb^\times.
$$ 
}
\end{defn}

\begin{prop}[\cite{cardy-alg}]
\label{prop:modinv-dim}
Let $A_\cl$ be a haploid commutative symmetric Frobenius algebra in $\CcC$.
If $A$ is modular invariant, then $A_\cl$ is also separable.
\end{prop}

\begin{prop}[\cite{cardy-alg}]  \label{prop:Aop-special}
Let $(A_\op|A_\cl, \iota_\co)$ be a Cardy $\Cc|\Cc^{(2)}$-algebra.
If $A_\cl$ is simple and $\dim A_\op \neq 0$,  
then $A_\op$ is simple and separable. 
\end{prop}

Now we give the reconstruction theorem for the Cardy $\Cc|\Cc^{(2)}$-algebra.
\begin{thm}[\cite{cardy-alg}]  \label{thm:reconst}
Let $A$ be a separable symmetric Frobenius $\Cc$-algebra. The triple $(A|Z(A), e)$, where $e: Z(A) \hookrightarrow R(A)$ is the canonical embedding, is a Cardy $\Cc|\Cc^{(2)}$-algebra. 
\end{thm}

\begin{rema} \label{rema:FRS}  {\rm
That an open-closed rational CFT can be constructed from a separable symmetric Frobenius algebra in $\Cc_V$ was first obtained by Fuchs, Runkel and Schweigert in a series of papers \cite{tft1}\cite{tft3} \cite{tft4}\cite{tft5},  in which the rational CFT is studied as a holographic boundary of a 3-dimensional topological field theory. In particular, they start from a separable symmetric Frobenius algebra in a modular tensor category $\Cc$ and use 3-dimensional TFT techniques to construct a so-called {\it solution of sewing constraints} \cite{tft5}, the notion of which can be proved to be equivalent to the notion of Cardy $\Cc|\Cc^{(2)}$-algebra \cite{cardy-alg-2}. I will refer to their approach to rational CFT as the FRS-framework.
}
\end{rema}

\noindent {\bf Examples}:  $X \otimes X^\vee$ for any $X\in \Cc$ are separable symmetric Frobenius algebras in $\Cc$. 
\bnu
\item Let $A=\one$, then $Z(\one) \cong \oplus_{i\in I} U_i^\vee \times U_i$.  
\item Let $A= X \otimes X^\vee$, then $Z(A) \cong Z(\one)$ as Frobenius algebras. 
\enu
$A_\cl=Z(\one)$ case is also called the Cardy case. We will show later that all $V$-invariant D-branes associated to $Z(\one)$ are given by $X \otimes X^\vee$. 
This is an indication of Morita equivalence which will be discussed later. 

\begin{rema}  \label{rema:cl-cft-mod=dbrane}  {\rm
In Cardy case ($A_\cl=Z(\one)$), all $V$-invariant D-branes are parametrized by objects $X$ in $\Cc$ and obtained by constructing the internal hom $[X, X]=X \otimes X^\vee$. The category of $V$-invariant D-branes has objects $X$ in $\Cc$ as objects and $[X, Y]=Y\otimes X^\vee$ as morphisms. Moreover, the functor $Y \mapsto [X, Y]$ is an equivalence between $\Cc$ and the category of $[X,X]$-modules \cite{ostrik}. Actually, one can show that the category of $Z(\one)$-modules (internal to $\CcC$) is monoidally equivalent to $\Cc$ \cite{eno2}. Therefore, in this case, all D-branes are also parametrized by the modules over the closed CFT. This is not true for cases other than Cardy. But to some extent, other cases in rational CFT can all be viewed as certain twists of the Cardy case by some internal symmetries \cite{ostrik}. 

On the other hand, we can also examine the chiral modules over $Z(\one)$ (recall Remark \ref{def:chiral-module}). A chiral $Z(\one)$-module that respects the chiral symmetry given by a rational VOA $V$ is nothing but a module in $\Cc$ over the algebra $T(Z(\one))=\oplus_{i\in I} U_i^\vee \otimes U_i$. It is easy to show that the category of chiral modules over $Z(\one)$ is equivalent to $\Cc^{\oplus |I|}$. Therefore, the category of boundary conditions is given by an indecomposable component of the category of chiral modules over $A_\cl$. 

A generalization of this result to cases other than Cardy is given in Theorem \ref{thm:open-exists}. For irrational CFTs, one can still use the notion of chiral module because it is defined at the level of vertex operators without using tensor category language. The existence of the OPE of vertex operators predicts that a certain weak version of tensor product should exist even in the irrational CFTs. I believe that $V_\op$ should be obtained from $M$ again by a certain internal-hom type of construction because the internal hom is nothing but the right adjoint functor of the tensor product. This should be the correct picture for irrational theories. 
}
\end{rema}

\bigskip
\section{Holographic principle and defects}  \label{sec:basic-SAG}

In Section \ref{sec:cl-cft}, we have shown that a closed CFT is stringy commutative. By our philosophy, it provides a new geometry SAG. In this section, we will use rational CFTs to see some basic properties of SAG. 

\subsection{Holographic Principle}  \label{sec:holography}

The Holographic Principle says that the information of the universe is contained in a part of it. In particular, such parts can have codimension higher than $1$ or even be a point. In SAG, a D-brane plays the role of a generalized point. So it is natural to ask if a single $V$-invariant D-brane determines the bulk theory uniquely. Indeed, we have the following uniqueness theorem: 
\begin{thm}[\cite{unique}\cite{cardy-alg}]  \label{thm:unique}
Let $(A|A_\cl,\iota_\co)$ be a Cardy $\Cc_V|\CcC$-algebra such that
$\dim A\neq 0$ and $A_\cl$ is simple. Then $A$ is separable and
$(A|A_\cl, \iota_\co) \cong (A|Z(A), e)$ as Cardy algebras.
\end{thm}

\begin{rema} {\rm
The above uniqueness result was first obtained by Fjelstad, Fuchs, Runkel and Schweigert \cite{unique} in the FRS-framework. In the context of Cardy $\Cc_V|\CcC$-algebra, it was proved in \cite{cardy-alg}. The equivalence of the two approaches is explained in \cite{cardy-alg-2}.
}
\end{rema}

It is unclear if the above uniqueness result still holds for general conformally invariant D-branes. But we expect that it still holds at least for certain nice conformally invariant D-branes.

\medskip
Conversely, for a closed CFT (or a bulk theory), we would like to know if there exists at least one D-brane and if it is unique. These questions are very important for us. The existence of D-branes means that the SAG associated to a given closed CFT is not empty, and the non-uniqueness means that the SAG contains more than one points. In the rational theories, we have the following result of the existence of $V$-invariant D-branes. 
\begin{thm}[\cite{cardy-alg}]\label{thm:open-exists}
If $A_\cl$ is a simple modular invariant commutative symmetric Frobenius algebra in $\CcC$, then there exist a simple 
separable symmetric Frobenius algebra $A$ in $\Cc$ and a morphism $\iota_\co: A_\cl \rightarrow R(A)$ such that  
\bnu
\item $A_\cl \cong Z(A)$ as Frobenius algebras;
\item $(A|A_\cl, \iota_\co)$ is a Cardy $\Cc_V|\CcC$-algebra;
\item $T(A_\cl) \cong \oplus_{\kappa\in J}\, 
M_{\kappa}^{\vee} \otimes_{A} M_{\kappa}$ as algebras, where $\{ M_{\kappa}\}_{\kappa \in J}$ is a set of representatives of the isomorphism classes of simple $A$-left modules. 
\enu
\end{thm}

\begin{rema}  \label{rema:non-cardy}  {\rm 
Notice that $V$-invariant D-branes can be parametrized by $\Cc_A$, the category of $A$-modules. On the other hand, $\Cc_{T(A_\cl)} \cong \Cc_A^{\oplus |J|}$.  Therefore, the category of boundary conditions is again an indecomposable component of the category of chiral modules over $A_\cl$. 
}
\end{rema}

In the Cardy case (Section \ref{sec:examples}), we see that D-branes are not unique. The ambiguity is controlled by the so-called Morita equivalence.  
\begin{defn}
Two algebras $A$ and $B$ in a tensor category $\Cc$ are called Morita equivalent if there are an $A$-$B$-bimodule $P$ and a $B$-$A$-bimodule $Q$ such that $P\otimes_B Q \cong A$ and $Q\otimes_A P \cong B$ as bimodules. 
\end{defn}
In the Cardy case, all $V$-invariant D-branes $X^\vee \otimes X$ are Morita equivalent to $\one$. Another example of Morita equivalence is given in Theorem \ref{thm:open-exists}, where all algebras $M_{\kappa}^{\vee} \otimes_{A} M_{\kappa}$ for $\kappa \in \mathcal{J}$ are Morita equivalent to $A$. Another equivalent way to define Morita equivalence is by the equivalence of the categories of modules over these two algebras. For an $A$-module $X$, $[X,X]$ is an algebra Morita equivalent to $A$ and the functor $Y \mapsto [Y, X], \forall \, Y\in \Cc_A$ is an equivalence between these two categories \cite{ostrik}. In general, we have the following result. 
\begin{thm}[\cite{cardy-alg}] \label{thm:morita}
If $(A_\op^{(i)}|A_\cl^{(i)}, \iota_\co^{(i)}), i=1,2$ are two Cardy $\Cc_V|\CcC$-algebras such that $A_\cl^{(i)}$ is simple and 
$\dim A_\op^{(i)}\neq 0$ for $i=1,2$, then $A_\cl^{(1)} \cong A_\cl^{(2)}$ as algebras if and only if $A_\op^{(1)}$ and $A_\op^{(2)}$ are Morita equivalent. 
\end{thm}

\begin{rema}  {\rm
The above theorem is equivalent to the statement that two simple separable symmetric Frobenius algebras $A$ and $B$ in $\Cc_V$ are Morita equivalent if and only if $Z(A) \cong Z(B)$ as Frobenius algebras. In other words, a single $V$-invariant D-brane determines uniquely the bulk theory; and the bulk theory determines uniquely the Morita class of the boundary theories.  
}
\end{rema}

Since a $V$-invariant D-brane is very symmetric and thus very scarce, we cannot see a continuum geometry here. I believe that by breaking the chiral symmetry down to conformal symmetry, one should be able to recover a rich and continuum geometry. For general conformally invariant D-branes, we need deal with irrational theories. The above strong result (Theorem \ref{thm:open-exists}) is perhaps not true for all conformally invariant D-branes. To study the general cases, it is useful to call two open-string VOAs {\it quasi-Morita equivalent} if they are both conformally invariant D-branes of the same closed CFT. It is important to understand this quasi-Morita equivalence in SAG. In some sense, such an equivalence can be viewed as the inverse statement of Holographic Principle which provides the foundation of SAG.

\subsection{Beauty and the Beast}   \label{sec:defect-duality}

We will briefly discuss a result of CFTs with defects that might be relevant to SAG. For more details of this subject, readers can consult \cite{ffrs-duality}\cite{runkel-suszek} and another contribution to this book \cite{dkr2}. 

\medskip
As we discussed in the introduction, we expect that spacetime are emergent from the structures of a CFT. For example, it is known to physicists that a space can emerge as the moduli space of D$0$-branes \cite{aspinwall}, which are certain 0-dimensional objects in the category of D-branes, or emerge in a certain limit of a family of CFTs \cite{konts-soib}\cite{soibelman}. What about time? Alain Connes observed that there is a God given embedding of the real number $\Rb$ in the group of outer automorphisms of a von Neumann algebra factor of type III \cite{connes}. He believes that it should be understood as a time evolution. It appeared later as Connes-Rovelli's Thermo Time Hypothesis \cite{time}. This deep observation seems to suggest that the emergence of time in SAG might also be related to the automorphism group (or duality group\footnote{In physics, duality usually can be weaker than an automorphism. But in this work, by duality we mean a true automorphism.}) of a CFT.

It has been known in many QFTs that defects play mysterious roles in various dualities \cite{savit}. Their relation has been clarified in the framework of rational CFTs by J\"{u}rg Fr\"{o}hlich, J\"{u}rgen Fuchs, Ingo Runkel and Christoph Schweigert \cite{ffrs-duality}. In particular, they show explicitly how an invertible defect can produce a duality of the bulk CFT\footnote{They also studied more general defects (called group-like) defects and their relation to non-invertible dualities.}. A $V$-invariant defect in a CFT over a rational VOA $V$ is a defect line separating two bulk phases and respecting the chiral algebra $V$ \cite{ffrs-duality}. If two bulk phases are determined by their boundary conditions, i.e. two simple separable symmetric Frobenius algebras $A$ and $B$ in $\Cc_V$, then a $V$-invariant defect is given by an $A$-$B$-bimodule $M$. Such a bimodule is invertible if there is another $A$-$B$-bimodule $N$ such that $M\otimes_B N \cong A$ and $N\otimes_A M\cong B$ as bimodules. The equivalence classes of invertible $A$-$A$-defects form a group called the Picard group denoted by $\mathrm{Pic}(A)$. A $V$-invariant automorphism of a bulk CFT is defined to be an automorphism of the bulk algebra $A_\cl=Z(A)$ such that it is the identity map on $V\otimes V$. We denote such $V$-invariant automorphism group by $\mathrm{Aut}(A_\cl)$.  Then we have the following exact correspondence between Beauty (the automorphism group) and the Beast (the invertible defects).

\begin{thm}[\cite{dkr}]  \label{thm:dkr}
For a simple separable symmetric Frobenius algebra $A$ in a modular tensor category $\Cc$, we have 
$\mathrm{Aut}(Z(A)) \cong \mathrm{Pic} (A)$.
\end{thm}

\begin{rema}  {\rm
The correspondence between duality and defects is not an isolated phenomenon. For example, a categorified version of Theorem \ref{thm:dkr} is proved in \cite{enom}\footnote{This result was obtained independently by Kitaev and myself \cite{kitaev-kong} and was announced in May 2009 in a conference on TQFT held at Northwestern University. In \cite{kitaev-kong} we identify the elements in $\mathrm{Pic}(\Cc)$ as the physical invertible defects of a lattice model and the simple objects in $Z(\Cc)$ as anyonic excitations in the bulk phase.}. More precisely, for a finite fusion category $\Cc$ we also have $\mathrm{Aut}(Z(\Cc))\cong \mathrm{Pic}(\Cc)$, where $Z(\Cc)$ is the monoidal center of $\Cc$. This can be viewed as a precise correspondence between duality and defects in Turaev-Viro TQFT (or Levin-Wen models) \cite{kitaev-kong}. The relation between duality and defects is known in many other QFTs \cite{kapustin-witten}\cite{kapustin-tikhonov}\cite{dgg}. I believe that under certain finite conditions, such as Hopkins-Lurie's fully dualizable condition \cite{lurie-tft}, an exact correspondence between dualities and invertible defects is true in many QFTs. 
}
\end{rema}

On the practical side, the invertible defects are much easier to compute than the dualities of a bulk theory. On the philosophical side, it suggests that the automorphism group of a QFT can emerge from the moduli of certain invertible defects. I believe that we can recover continuum automorphism groups if we study conformally invariant defects. It suggests that the emergence of time might be related to invertible defects. Notice that we have only discussed the automorphisms of the closed CFT. How about the automorphism groups of an open CFT? It is very possible that each D-brane carries its own time evolution. Moreover, it was proposed by Connes and Rovelli \cite{time} in their Thermo Time Hypothesis that a time involution depends on the choice of a thermo-state. So it is desirable to define thermo-states or KMS states in the context of CFT. We will return to these questions in our future studies.

\medskip
Because the Holographic Principle is so important in our program, it is worthwhile to explore its properties further. Since a boundary theory often determines the bulk theory by taking the center, we would like to ask if this process is functorial. By considering the groupoid of those algebras $A$ described in the Theorem \ref{thm:dkr} and invertible bimodules,  and the groupoid of commutative algebras in $Z(\Cc)$ and isomorphisms of algebras, we obtain a groupoid version of Theorem \ref{thm:dkr} \cite{dkr}. Namely, it is an equivalence of groupoid. It turns out that this equivalence of groupoid can be extended to a richly structured functor which take an algebra to its center. Namely, the notion of the center is functorial. This suggests that the Holographic Principle is also functorial. More interestingly, this functoriality demands
defects of all codimensions to appear. As a purely algebraic statement, this functoriality properly encodes all the information of defects without attributing it to a geometric bordism category. This seems to suggest a purely algebraic way to understand extended TQFT \cite{lurie-tft} (see also \cite{bergner} in this book). See \cite{dkr2} in this book for more details on the functoriality of the center and its relation to TFTs with defects.

\bigskip
\section{Conclusions and outlook}  \label{sec:outlooks}

In this work, we have outlined a dream of a new algebraic geometry (SAG) based on 2-dimensional conformal field theories. We also reviewed some recent progress on the mathematical foundations of rational open-closed CFTs. Our understanding on $V$-invariant D-branes in a rational CFT is quite satisfying. But the real challenge lies in how to understand conformally invariant D-branes in concrete examples in order to see a rich geometry. In this paper, we did not touch upon open-closed superconformal field theory which is the real hero behind many miracles in geometry.  But its mathematical foundation is still lacking. We will return to these issues in the future publications.

\medskip
Actually, the philosophy of SAG itself can have useful applications in other fields. Many structures that have appeared in CFT also appeared in many different contexts. For example, the Holographic Principle appears in the context of $E_n$-algebra or $E_n$-category as the so-called generalized Deligne conjecture \cite{nlab-deligne}; defects reincarnate in algebraic geometry as Fourier-Mukai transformations, etc. Therefore, questions asked here can also be transported to other fields, and vice versa.  Let us discuss two concrete cases below.

\begin{itemize}

\item String topology was invented by Chas and Sullivan \cite{string-top}. It can be viewed as an open-closed homological CFT \cite{godin}.  Let $M$ be a simply connected and closed manifold. The closed string topology is an algebraic structure on $H_\ast(LM)$, where $LM$ denotes the loop space. This is the closed algebra in this case. By our philosophy of SAG, we would like to ask what lies in the spectrum. Namely, what are the compatible open algebras? Can they recover all the ``points" in $M$? Let $N$ be a submanifold of $M$. The chain space $C_\ast(\mathcal{P}_{N, N})$ of the path space between $N$ and $N$ is a differential graded algebra. On the level of homology, this algebra structure is the string topology product introduced by Sullivan \cite{sullivan}. The dg-algebra $C_\ast(\mathcal{P}_{N, N})$ is an open algebra in this case. The notion of center in this case is given by Hochschild cohomology. Then our question on the spectrum of the closed string topology is related to the following question asked by Blumberg, Cohen and Teleman \cite{bct}: For what submanifolds $N$ the relation: 
$$
HH^\ast(C_\ast(\mathcal{P}_{N, N}), C_\ast(\mathcal{P}_{N, N})) \cong H_\ast (LM)
$$
holds? As pointed out in \cite{bct}, the answer is affirmative for a point in $M$ and $N=M$ by works in the 1980s.  They went on to give more answers to their question including the inclusion $N \hookrightarrow M$ being null homotopic or being the inclusion of the fiber of a fibration $M \to B$. One can also ask many other obvious questions. For example, can one give a geometric meaning (as some kind of a boundary condition) of a certain module over the closed algebra $H_\ast (LM)$ so that the open algebra $C_\ast(\mathcal{P}_{N, N})$ arises as the internal hom of this module? 

\medskip
\item A crucial result in the local geometric Langlands correspondence (GLC) says that the center of the vertex algebra $V_{-h^\vee}(\mathfrak{g})$ associated to a simple Lie algebra $\mathfrak{g}$ at the critical level $k=-h^\vee$ is isomorphic to the function algebra over the space of ${}^L\mathfrak{g}$-opers on the formal disk \cite[Theorem 9]{edward-frenkel}. Actually, the critical level is exactly the case of non-CFT because the Sugawara construction of Virasoro algebra fails exactly in this case. The vertex operator algebra $V_k(\mathfrak{g})$ for $k\neq -h^\vee$ can be viewed as a D-brane in the Cardy case. Its center gives the closed CFT $V_\cl$. Moreover, by SAG, the spectrum of this closed CFT is just the category of D-branes. Therefore, we expect what appears on the one side of GLC for noncritical levels is the category of D-branes of the closed CFT $V_\cl$ associated to $\mathfrak{g}$. It is very possible that the other side of GLC is also given by the category of D-branes of a closed CFT associated to ${}^L\mathfrak{g}$. This is supported by the intuition from physics. The physics origin of GLC is a conjectured duality, which is derived from Montonen-Olive duality by topological twisting, between two 4-dimensional topological gauge theories with gauge group $G$ and ${}^LG$ respectively \cite{kapustin-witten}. One consequence of this duality is a conjectured equivalence between two 3-categories of boundary conditions of these two topological gauge theories \cite{kapustin}. 

\end{itemize}

\medskip
Although the speculation of a new geometry SAG discussed in this paper is still too naive and premature, it indeed motivated some of my own works on CFTs. I hope that this naive picture can inspire more serious works in this direction. Another reason for me to write about it is to respond to some opinions I heard in various occasions that 2-dimensional CFTs are well understood, and there are not many interesting things left to do. I hope that this paper can convince some of my readers that there are still a lot of interesting questions in 2-dimensional CFT waiting to be studied.


\bibliographystyle{amsalpha}

\end{document}